\documentclass[10pt,a4paper]{article}
\usepackage[utf8]{inputenc}
\usepackage{amssymb,amsmath,dsfont}
\usepackage{graphicx}
\usepackage{wrapfig}
\usepackage{theorem} 
\usepackage[margin=3cm]{geometry}
\usepackage[dvipsnames,usenames]{color}
\usepackage{epstopdf}
\usepackage[font=small,labelfont=bf]{caption}
\usepackage{hyperref}
\usepackage{comment}
\usepackage[pagewise]{lineno}

\newtheorem{Theorem}{Theorem}[section]

\newtheorem{Proposition}{Proposition}[section]
\newtheorem{Lemma}{Lemma}[section]

\newtheorem{Remark}{Remark}[section]

{\theorembodyfont{\rmfamily}\newtheorem{test}{Test}}

\title{Reconstruction of degeneracy region and power for parabolic equations and systems}
\author{Piermarco Cannarsa\thanks{Department of Mathematics, University of Rome Tor Vergata, Italy (e-mail: {\tt cannarsa@mat.uniroma2.it}).}\, , Veronica Danesi\thanks{Department of Mathematics, University of Rome Tor Vergata, Italy (Corresponding author, e-mail: {\tt danesi@mat.uniroma2.it}).}\, , Anna Doubova\thanks{Dpto. EDAN and IMUS, Universidad de Sevilla, Spain (e-mail: {\tt doubova@us.es}).} }
\date{\today}

\begin{document}

\maketitle

\begin{abstract}
    We address the inverse problem of recovering a degeneracy point within the diffusion coefficient of a one-dimensional complex parabolic equation by observing the normal derivative at one point of the boundary.
    The strongly degenerate case is analyzed.
    In particular, we derive sufficient conditions on the initial data that guarantee the stability and uniqueness of the solution obtained from a one-point measurement. Moreover, we present more general uniqueness theorems, which also cover the identification of the initial data, the coefficient of the zero order term and the degeneracy power, using measurements taken over time.
    Our method is based on a careful analysis of the spectral problem and relies on an explicit form of the solution in terms of Bessel functions.
    Our investigation also covers the case of real 1-D degenerate parabolic systems of equations coupled with a specific structure.
    Theoretical results are also supported by numerical simulations.
\end{abstract}

\begingroup
\renewcommand{\thefootnote}{}
\addtocounter{footnote}{-2}
\endgroup


\section{Introduction}

The aim of this paper is to investigate the inverse problem of reconstructing an interior degeneracy point $a\in (0,1)$ for the following degenerate parabolic complex equation
\begin{equation}\label{eq.deg}
\left\{\begin{array}{ll}
\partial_t w- \partial_x(|x-a|^\theta \partial_x w) -cw= 0, & (x,t)\in(0,1)\times(0,T),  \\[1mm]
w(0,t) = 0, \quad w(1,t) = 0, &  t\in (0,T),  \\[1mm]
w(x,0) = w_0(x),  & x\in(0,1),
\end{array}\right.
\end{equation}
where $T>0$, $\theta \in [1,2)$, $w_0(x)=u_0(x)+i v_0(x)\neq 0$ with $u_0,v_0$ real-valued functions, $c=\alpha + i \beta$, with $\alpha,\beta\in \mathbb{R}$, are given. 
Specifically, we consider the strongly degenerate case with $1\leq\theta<2$ in the diffusion coefficient.

Our goal is to determine or approximate the degeneracy point $a\in(0,1)$ from suitable measurements of the solution.
The unknown degeneracy point being inside the domain, a natural extra observation of the solution to the above problem is the normal derivative $\partial_x w(x,t)$ at the boundary.
In particular, this simplified model describes heat diffusion in a body with a conductivity failure. 
The diffusion coefficient is usually related to the structure of the material, the density, and other factors.
Thus, the degeneracy of this coefficient indicates the ability to resist heat transfer.
The objective is then to determine the unknown location of this degeneracy using suitable boundary heat flux data.

Degenerate parabolic equations have attracted increasing attention due to their significant theoretical implications and wide-ranging practical applications in fields such as climatology (see \cite{Diaz,Ji-Huang,Sellers}), financial mathematics (see \cite{Black-Scholes}), fluid dynamics (see \cite{Oleinik-etal}), and population genetics (see \cite{Ethier}).
Despite their theoretical and practical importance, the literature concerning inverse problems for degenerate parabolic PDEs is relatively new. 
Examples include the inverse source problem (see \cite{PMC-Mrt-Van-16,PMC-Tort-Yam-10,Deng-etal-15,Hussein-etal-20,Kamynin-18,Li-Li-21,Tort}), the recovery of the first-order coefficient (see \cite{Deng-Yang,Kamynin-20}), or various identification problems of degenerate diffusion coefficients, which also encompass the reconstruction of the power exponent (see \cite{PMC-Dou-Yam-21}).
These degenerate problems can also be divided into different classes according to the way of degeneration with respect to either spatial variables or to the time variable.
For instance, works addressing the reconstruction of a time-dependent degenerate diffusion coefficient, and similarly, the recovery of a time-dependent first-order coefficient, can be found in \cite{Huzyk-etal-23,Ivanchov-etal-06}.

For example, the inverse problem of reconstructing an interior degeneracy point for the real case was considered in \cite{PMC-Dou-Yam-24}, where the authors analyzed the strongly degenerate case for $\theta=1$.
Our goal is to generalize this result in several directions. 
One direction involves extending it to systems of two real degenerate coupled parabolic equations. 
Specifically, in this work, we consider a coupling with a particular structure that allows us to reformulate the problem as a complex degenerate parabolic equation.
Another direction of generalization concerns extending the result to cases where the degeneracy has an exponent different from 1. 
In particular, we will analyze the strongly degenerate case with $\theta\in[1,2)$.
Here, the analysis is technically more complicated and our results generalize those in \cite{PMC-Dou-Yam-24}.
Instead, the weakly degenerate case with $\theta\in[0,1)$ is still an open problem. 
In fact, we cannot analyze the spectral behavior independently to the left and right of the degeneracy. 
This distinction prevents a spectral analysis analogous to what is possible in the strongly degenerate situation, with the exception of specific configurations of the degeneracy point.
Furthermore, from an applicative point of view, as in climatology (see \cite{Kaper}) or financial mathematics (see \cite{Baldi}, Chapter 13, Section 6), the most interesting examples fall into the class of strongly degenerate problems.
It should also be noted that the restriction $\theta<2$ is due to the spectral technique implemented, specifically to the use of Bessel functions.

The last extension regards the uniqueness and reconstruction for the inverse problems where also the degeneracy power, the coefficient of the zero-order term, and the initial data are unknown.
In this case, our result for the degeneracy power also differs from the one in \cite{PMC-Dou-Yam-21}, in terms of the measurements considered and for the extension to systems. 
In particular, in \cite{PMC-Dou-Yam-21} the authors take into account measurements of the time derivative and space gradient for a fixed time, but distributed on the entire spatial domain.
Our approach considers a more natural observation, represented by the normal derivative at one boundary distributed over a time interval.

Inverse problems are classified as ill-posed in the Hadamard sense.
This means that their solution may not exist, may be non-unique, and/or may be highly sensitive to small errors in the provided data, leading to significant inaccuracies in the computed solutions.
The main issues concerning our interior degeneracy reconstruction problem are uniqueness, stability, and numerical approximation of the solution.
For this last aim, we will reformulate our inverse problem as an optimization problem (see Section~\ref{sec.numerics}). 
This approach is a standard technique in inverse problems for reconstructing unknown data, and similar methods have been used in previous studies for other problems (see \cite{Apraiz-etal,PMC-Dou-Yam-21,Dou-etal-15, Lavrentiev-etal, Samarskii-etal, Vogel}).

The paper is organized as follows. In Section 2, we introduce the functional setting and establish the well-posedness of the corresponding direct problem. 
In Section 3, we analyze the eigenvalue problem.
In Section 4, we provide an expression of the normal derivative computed using Bessel
functions. 
In Section 5, we establish a Lipschitz stability result with one-point measurements. 
Section 6 is devoted to general uniqueness results for distributed measurements over a time interval, assuming that the initial data, the coefficient $c$ and the degeneracy power are
also unknown.
Section 7 concerns the application of the previous results to real systems of 1-D coupled degenerate parabolic equations.
Section 8 concludes with numerical experiments related to the inverse problem under consideration.

\section{Functional setting and well-posedness}

In this Section, we introduce the appropriate weighted energy spaces in which the problem can be set, depending on the value of the parameter $\theta$.
Moreover, the well-posedness of the direct problem will be stated.

Consider $X = L^2 (0,1;\mathbb{C}) $
endowed with the scalar product $\langle f , g \rangle = \int _{0} ^1 f(x)\Bar{g}(x)\; dx$, $\forall\; f,g \in   X$.
We define
\begin{equation*}
\begin{aligned}
H^1_{\theta} (0,1;\mathbb{C}):=
 &\Bigl \{ w \in X \ \mid \  w \text{ locally absolutely continuous in } (a,1] \text{ and }
 \text{in } [0,a),\\ 
 & \qquad\int _{0} ^1 \vert x-a \vert ^\theta \vert w^\prime(x)\vert ^2  \, dx < \infty
\text{ and } w(0)=0 =w(1) \Bigr \} , \quad 1\leq\theta<2 ,
\end{aligned}
\end{equation*}
that is endowed with the natural scalar product
$$ (f,g) = 
\int _{0} ^1  \bigl( \vert x-a \vert ^\theta f'(x) \Bar{g}'(x) + f(x) \Bar{g}(x) \bigr) \, dx\ , \quad  \forall f,g \in  H^1_{\theta} (0,1;\mathbb{C}).$$
Next, consider 
\begin{equation*}
H^2_{\theta} (0,1;\mathbb{C}):=   \Bigl \{ w \in H^1_{\theta }(0,1;\mathbb{C})  \ \mid \int _{0} ^1 \vert (\vert x-a \vert ^\theta w'(x) ) ' \vert ^2 \, dx < \infty  \Bigr \} , \quad 1\leq\theta<2,
\end{equation*}
and the operator $\mathbb{A}:D(\mathbb{A})\subset X \to X$ will be defined by
$D(\mathbb{A}) :=  H^2_{\theta} (0,1;\mathbb{C})$ 
and
\begin{equation*}
\mathbb{A}w:=Au+iAv\, \quad \forall w=u+iv \in D(\mathbb{A})\quad \text{and}\quad u,v \;\;\mathbb{R}-\text{valued functions}, 
\end{equation*}
with
\begin{equation*}
    A:=\partial_x(\vert x-a \vert ^\theta \partial_x ),\qquad
    D(A) :=  H^2_{\theta} (0,1)=\Bigl \{ u \in H^1_{\theta }(0,1)  \ \mid \int _{0} ^1 \vert (\vert x-a \vert ^\theta u'(x) ) ' \vert ^2 \, dx < \infty  \Bigr \},
\end{equation*}
and
\begin{equation*}
\begin{aligned}
    H^1_{\theta} (0,1):=
 &\Bigl \{ u \in L^2 (0,1) \ \mid \  u \text{ locally absolutely continuous in } (a,1] \text{ and } 
\text{in } [0,a),\\ 
&\qquad\int _{0} ^1 \vert x-a \vert ^\theta \vert u^\prime(x)\vert ^2  \, dx < \infty
\text{ and } u(0)=0 =u(1) \Bigr \} , \quad 1\leq\theta<2 .
\end{aligned}
\end{equation*}
Then, the following results hold:
\begin{Proposition}
\label{Prop-A-s}
Given $\theta \in [1,2)$, we have:

a) $H^1_{\theta} (0,1;\mathbb{C})$ is a Hilbert space.

b) $\mathbb{A}: D(\mathbb{A}) \subset X \to X$ is a dissipative self-adjoint operator with dense domain.

Therefore, $\mathbb{A}$ is the infinitesimal generator of an analytic semigroup of contractions $e^{t\mathbb{A}}$ on $X$ and $t\mapsto w(\cdot,t)$ is an analytic map for all $t>0$.
\end{Proposition}
\noindent
 \textbf{Proof of Proposition~\ref{Prop-A-s}:}
       The proof of a) and b) is similar to that of the real case in \cite{PMC-Frr-Mrt-19}. 
       Analyticity follows from a well-known result (see \cite{Tanabe}).
 \hfill$\blacksquare$
\vskip 4pt
 
Given an initial condition $w_0 \in  X$, the problem \eqref{eq.deg} can be recast in the abstract form
\begin{equation}
\label{abstract-pb}
    \begin{cases}
 w'(t)=(\mathbb{A}+c\mathbb{I})w(t) & t\ge 0,
 \vspace{.1cm}
 \\
 w(0)=w_0.
\end{cases}
\end{equation}
The function $w \in \mathcal C ^0 ([0,T]; X)  \cap L^2(0,T; H^1_{\theta} (0,1;\mathbb{C}))$,
given by the formula
$$ w(\cdot ,t) = e^{t(\mathbb{A}+c\mathbb{I})} w_0=e^{(\alpha+i\beta)t}\left(e^{tA}u_0+ie^{tA}v_0 \right),
$$ 
is the solution of \eqref{abstract-pb} in the sense of semigroup theory. 
We say that a function 
$$w \in
\mathcal C ^0 ([0,T]; H^1_{\theta} (0,1;\mathbb{C}))  \cap   H^1(0,T; X) \cap    L^2 (0,T; D(\mathbb{A})) $$
is a \textit{strict} solution of \eqref{abstract-pb} if $w$ satisfies $\partial_t w - \partial_x(\vert x-a \vert ^\theta \partial_x w) -cw= 0$ almost everywhere
in $(0,1)\times (0,T)$, and the initial and boundary conditions for all $t\in [0,T]$ and all $x\in [0,1]$.  
Moreover, it is possible to prove the existence and uniqueness of the strict solution.
In particular, the following result holds true.

\begin{Proposition}
\label{prop-wp}
If $w_0 \in H^1 _{\theta}(0,1;\mathbb{C})$, then the mild solution of \eqref{abstract-pb} is the unique strict solution of \eqref{abstract-pb}.
\end{Proposition}

\noindent
 \textbf{Proof of Proposition~\ref{prop-wp}:}
        The proof is analogous to that in the real case (see, for instance, \cite{PMC-Mrt-Van-18}).
  \hfill$\blacksquare$
\vskip 5pt

\begin{Remark}\label{Neumann-cond}
We also observe that, for $w\in H^2_\theta(0,1;\mathbb{C})$ and $\theta\in[1,2) $, we have $\vert x-a\vert^\theta \partial_x w|_{x=a}=0$. 
In fact, if $\vert x-a\vert^\theta \partial_x w(x) \to L$ when $x\to a$, then $\vert x-a\vert^\theta|\partial_x w(x) |^2 \sim L^2/\vert x-a\vert^\theta$ and therefore $L=0$, otherwise $w\notin H^1_\theta(0,1;\mathbb{C})$.

As a consequence, the strongly degenerate problem can be decoupled into two completely distinct sub-problems. 
More specifically, in the strongly degenerate case, the two problems in $(0,a)$ and $(a,1)$ can be analyzed separately, taking into account the Neumann boundary condition in $x=a$.
\end{Remark}

\section{The eigenvalue problem}  

The analysis of the spectral problem associated to \eqref{eq.deg} will be essential for our purposes.
The eigenvalues and associated eigenfunctions of the degenerate diffusion operator $-A$ are the nontrivial solutions $(\lambda,\phi)\in\mathbb{R}\times D(A)$ of
\begin{equation}
\label{eigen-pbm}
\begin{cases}
 - (\vert x-a \vert ^\theta \phi '(x))' =\lambda \phi (x),  \qquad x\in (0,1),\\
\phi(0)=0 =\phi(1),
\end{cases}
\end{equation}
which can be expressed in terms of Bessel functions of the first kind (see \cite{Kamke}). 

\medskip

For $\theta \in [1,2)$, let
\begin{equation*}
\nu _\theta := \frac{ \vert \theta -1 \vert }{2-\theta} = \frac{ \theta -1 }{2-\theta}, 
\qquad k_\theta:= \frac{2-\theta}{2}.
\end{equation*}
Given $\nu_\theta$, we denote by $J_{\nu_\theta}$ the Bessel function of the first kind and of order $\nu_\theta$ given by
\begin{equation}\label{eq.besselJnu}
J_{\nu_\theta}(z) = \displaystyle\sum_{k=0}^\infty \dfrac{(-1)^k}{k! \ \Gamma(k+\nu_\theta+1)}\Big( \frac{z}{2}\Big)^{2k+\nu_\theta}, \quad z\ge 0,
 \end{equation}
where $\Gamma$ is the Gamma function (see \cite{Watson}).
Moreover, let us denote by $j_{\nu_\theta,1}< j_{\nu_\theta,2} < \ldots < j_{\nu_\theta,n} <\dots$ the positive zeros of $J_{\nu_\theta}$.

When $\theta \in [1,2)$, we have the following description of the spectrum of the associated operator:

\begin{Proposition}
\label{prop-eigenf}
The admissible eigenvalues $\lambda$ for problem \eqref{eigen-pbm} are given by 
\begin{equation*}
\forall n \geq 1, \qquad \lambda_{ n}^{(r)}(a) = k_\theta ^2 \frac{j_{\nu _\theta ,n}^2}{(1-a)^{2k_\theta}} \qquad \text{or} \qquad  \lambda_{ n}^{(l)}(a) = k_\theta ^2 \frac{j_{\nu _\theta ,n}^2}{a^{2k_\theta}}.
\end{equation*}
An orthonormal basis in  $L^2(0,1)$ is given by the following eigenfunctions
\begin{equation*}
\label{*fp-d}
\tilde \phi ^{(r)} _{\theta, n}(x) :=
 \begin{cases}
 0 \quad & \text{ if } x \in (0,a), \\
\dfrac{\sqrt{2 k_\theta }}{ \sqrt{1-a} \,\vert J'_{\nu_\theta} (j_{\nu_\theta,n} ) \vert} \left(\dfrac{x-a}{1-a}\right)^{\frac{1-\theta}{2}} J_{\nu _\theta} \left(j_{\nu_\theta,n} \left(\dfrac{x-a}{1-a}\right) ^{k_\theta}\right) \quad & \text{ if } x \in (a,1), 
\end{cases} 
\end{equation*}
and
\begin{equation*}
\label{*fp-g}
\tilde  \phi ^{(l)} _{\theta, n}(x) :=
 \begin{cases}

\dfrac{\sqrt{2 k_\theta }}{\sqrt{a}\, \vert J'_{\nu_\theta} (j_{\nu_\theta,n} ) \vert} \left| \dfrac{x-a}{a} \right| ^{\frac{1-\theta}{2}} J_{\nu _\theta} \left(j_{\nu_\theta,n} \left| \dfrac{x-a}{a} \right|  ^{k_\theta}\right)\quad & \text{ if } x \in (0,a), \\ 
0 \quad & \text{ if } x \in (a,1).
\end{cases}
\end{equation*}

\end{Proposition}

\noindent
 \textbf{Proof of Proposition~\ref{prop-eigenf}:}
The eigenvalue problem \eqref{eigen-pbm} can be split into the following two sub-problems
\begin{equation*}
\label{eigen-pbm-right}
\begin{cases}
 - (\vert x-a \vert ^\theta \phi '(x))' =\lambda \phi (x),  \qquad x\in (a,1),\\
\phi(1)=0,
\end{cases}
\end{equation*}
and
\begin{equation*}
\label{eigen-pbm-left}
\begin{cases}
 - (\vert x-a \vert ^\theta \phi '(x))' =\lambda \phi (x),  \qquad x\in (0,a),\\
\phi(0)=0,
\end{cases}
\end{equation*}
which can be transformed into the two following sub-problems 
\begin{equation*}\label{pb.autofunc-right}
 \begin{cases}
 - (y^\theta \varphi'(y))' = \lambda (1-a)^{2-\theta} \varphi(y),  \quad y\in (0,1),
 \\[2mm]
 \varphi(1) = 0,
 \end{cases}
 \end{equation*}
 and
 \begin{equation*}\label{pb.autofunc-left}
 \begin{cases}
 - (\vert y\vert ^\theta \varphi'(y))' = \lambda a^{2-\theta} \varphi(y),  \quad y\in (-1,0),
 \\[2mm]
 \varphi(-1) = 0,
 \end{cases}
 \end{equation*}
 by means of the two coordinate transformations $y = \dfrac{x-a}{1-a}$, with $\varphi(y)=\phi(a+(1-a)y)$, and $y = \dfrac{x-a}{a}$, with $\varphi(y)=\phi(a+ay)$, respectively.
 \vskip 2pt
 The first eigenvalue sub-problem can be rewritten as a differential Bessel's equation of order $\nu_\theta=\frac{\theta-1}{2-\theta}$
 \begin{equation}
\label{pb-bessel}
\begin{cases}
z^2 \Psi ''(z) + z \Psi  '(z) + (z^2 - \nu_\theta^2) \Psi (z) = 0, \quad z\in \left(0, \frac{2 }{2-\theta}\sqrt{\lambda}(1-a)^{\frac{2-\theta}{2}} \right), \\
\Psi \Bigl(\frac{2}{2-\theta} \sqrt{\lambda} (1-a)^{\frac{2-\theta}{2}}\Bigr) = 0, 
\end{cases}
\end{equation}
 by setting
    $\varphi(y):= y^{\frac{1-\theta}{2}}\Psi\left(\frac{2}{2-\theta}\sqrt{\lambda}((1-a)y)^{\frac{2-\theta}{2}} \right)$.
The second one leads to the Bessel's equation
\begin{equation*}
\begin{cases}
z^2 \Psi ''(z) + z \Psi  '(z) + (z^2 - \nu_\theta^2) \Psi (z) = 0, \quad z\in \left(0, \frac{2 }{2-\theta}\sqrt{\lambda}a^{\frac{2-\theta}{2}} \right), \\
\Psi \Bigl(\frac{2}{2-\theta} \sqrt{\lambda} a^{\frac{2-\theta}{2}}\Bigr) = 0, 
\end{cases}
\end{equation*}
 by setting
    $\varphi(y):= |y|^{\frac{1-\theta}{2}}\Psi\left(\frac{2}{2-\theta}\sqrt{\lambda}(a|y|)^{\frac{2-\theta}{2}} \right)$.
The proof follows from the result in \cite{PMC-Frr-Mrt-19}, taking into account the changes of coordinates $y=\dfrac{x-a}{1-a}$ and $y=\dfrac{x-a}{a}$, respectively, and with a suitable modification of the eigenvalues, determined by means of the boundary condition.
  \hfill$\blacksquare$

  \begin{Remark}
      Observe that in case $\theta=1$, the computation for the eigenfunctions in \cite{PMC-Dou-Yam-24} does not contain the factor $\dfrac{1}{\sqrt{1-a}}$. This extra term obviously has to be inserted in the computations. However, as we shall see in our case, no changes in the result will be necessary.
  \end{Remark}

 We now recall some properties of the Bessel functions that will be used later.
\begin{Lemma}[Properties of Bessel functions]\label{lemma.Bessel}
Let $J_\nu(z)$, with $\nu\in\mathbb{R}$, be the Bessel functions of order $\nu$ and of the first kind (given by \eqref{eq.besselJnu}) and let us denote by $\{j_{\nu,n}\}_{n\ge 1}$ the sequence of increasing positive zeros of $J_\nu$, i.e. $J_\nu(j_{\nu,n})=0$, with $0<j_{\nu,1}<j_{\nu,2}<\cdots$. 

Then, the following properties hold: 
\begin{enumerate}
\item[\textit{a)}] $ \dfrac{d}{dz}\big (z^\nu J_\nu(z)\big) =  z^\nu J_{\nu-1}(z)$;
\item[\textit{b)}] $J_{\nu-1}(z) + J_{\nu+1}(z) =\dfrac{2\nu}{z} J_\nu(z)$ and $zJ_{\nu}'(z)-\nu J_{\nu}(z)=-zJ_{\nu+1}(z)$;  
\item[\textit{c)}] $\vert J_\nu(z)\vert\leq 1$ for $\nu\geq 0$;
\item[\textit{d)}]  $\displaystyle\int_0^{j_{\nu,n}} s^{\nu+1} J_\nu(s) \, ds = j_{\nu,n}^{\nu+1} J_{\nu+1}(j_{\nu,n}) = -j_{\nu,n}^{\nu+1} J_\nu'(j_{\nu,n})$;
\item[\textit{f)}] $\forall\; \nu\in\left[0,\frac{1}{2}\right]$, $\forall\, n\geq 1$, 
$\pi\Big(n +\dfrac{\nu}{2}-\dfrac{1}{4}\Big) \le j_{\nu,n} \le \pi\Big(n+\dfrac{\nu}{4} -\dfrac{1}{8}\Big)$;
\item[\textit{g)}] $\forall\; \nu\geq\frac{1}{2}$, $\forall\, n\geq 1$, 
$\pi\Big(n +\dfrac{\nu}{4}-\dfrac{1}{8}\Big) \le j_{\nu,n} \le \pi\Big(n+\dfrac{\nu}{2} -\dfrac{1}{4}\Big)$.
\end{enumerate} 

\end{Lemma}

\noindent


\section{Computation of the normal derivative}\label{sec:normal-der}

In this section, we will perform the explicit computation of the normal derivative $\partial_x w^a(x,t)$ at the boundary, where $w^a(x,t)$ is the solution of \eqref{eq.deg}.
Notice that the function $t\mapsto\partial_x w^a(1,t)$ is analytic for all $t>0$, since $w^a$ is analytic for all $t>0$.
In the following, we will only consider the problem in the right interval $(a,1)$.
A similar analysis can also be performed in the left interval $(0,a)$, taking into account $\partial_x w^a(0,t)$.

\medskip

For the strongly degenerate case, we concentrate on the analysis of the sub-problem
\begin{equation}\label{eq.deg-right}
\left\{\begin{array}{ll}
\partial_t w- \partial_x((x-a)^\theta \partial_x w) -cw= 0, & (x,t)\in(a,1)\times(0,T),  \\[1mm]
 (x-a)^\theta \partial_x w|_{x=a}=0, \quad w(1,t) = 0, &  t\in (0,T),  \\[1mm]
w(x,0) = w_0(x),  & x\in(a,1),
\end{array}\right.
\end{equation}
where we have taken into account the Neumann boundary conditions in $x=a$ (see Remark \ref{Neumann-cond}).

With the aim of computing the normal derivative $\partial_x w^a(1,t)$, we introduce the following change of variables
 \begin{equation}\label{eq.change}
y = \dfrac{x-a}{1-a},\quad  x = a + (1-a)y,
 \end{equation}
obtaining
$w^a(x,t)=\tilde{w}\Big(\dfrac{x-a}{1-a}, t\Big)$, $x\in (a,1)$ and $\tilde{w}_0^a(y)  = w_0(a + (1-a)y)$, $y\in (0,1)$.

Therefore, $\tilde{w} = \tilde{w}(y,t)$ satisfies 
\begin{equation*}
\left\{\begin{array}{ll}
\partial_t \tilde{w}- \dfrac{1}{1-a}\partial_y(y^\theta\partial_y \tilde{w})-c\tilde{w} = 0, & (y,t)\in(0,1)\times(0,T),  \\[1mm]
y^\theta\partial_y \tilde{w}(y,t)\big|_{y=0} = 0, \quad \tilde{w}(1,t) = 0, &  t\in (0,T),  \\[1mm]
\tilde{w}(y,0) = \tilde{w}_0^a(y),  & y\in(0,1)
\end{array}\right.
\end{equation*}
and the normal derivative reads $\partial_x w^a(1,t) = \dfrac{1}{1-a}\partial_y \tilde{w}(1,t)$.
In addition, as we have shown in the previous section, the eigenvalue problem associated to the degenerate diffusion operator 
 \begin{equation}
 \label{spectral-pb-right}
 \begin{cases}
 - ((x-a)^\theta \phi'(x))' = \lambda \phi,  \quad x\in (a,1),
 \\[1mm]
 \phi(1) = 0,\quad   (x-a)^\theta\phi'(x)\big|_{x=a} = 0
 \end{cases}
 \end{equation}
can be transformed into a problem on the interval $(0,1)$ by the change of variables \eqref{eq.change}. 
We get that $\varphi(y) = \phi(a+(1-a)y)$ satisfies
  \begin{equation*}
 \begin{cases}
 - (y^\theta \varphi'(y))' = \lambda (1-a)^{2-\theta} \varphi,  \quad y\in (0,1),
 \\[1mm]
 \varphi(1) = 0,\quad   y^\theta\varphi'(y)\big|_{y=0} = 0.
 \end{cases}
 \end{equation*}
In particular, the eigenvalues are given by 
\begin{equation}
\label{eigenvalues}
    \lambda_{ n} = k_\theta ^2 \frac{j_{\nu _\theta ,n}^2}{(1-a)^{2k_\theta}}
\end{equation}
and the corresponding eigenfunctions can be written as
\begin{equation}\label{eq.autofun}
    \phi_n(x)=\dfrac{\sqrt{2 k_\theta }}{\sqrt{1-a}\, \vert J'_{\nu_\theta} (j_{\nu_\theta,n} ) \vert}G_{\nu_\theta,n}(x,a),
\end{equation}
where
\begin{equation*}
    G_{\nu_\theta,n}(x,a):=\left(\dfrac{x-a}{1-a}\right)^{\frac{1-\theta}{2}} J_{\nu _\theta} \left(j_{\nu_\theta,n} \left(\dfrac{x-a}{1-a}\right) ^{k_\theta}\right).
\end{equation*}
Introducing the notation $d_{\nu_\theta,n}:=J_{\nu_\theta}'(j_{\nu_\theta,n}) \left(j_{\nu_\theta,n}\right)^{\frac{1}{2k_\theta}}$, $h_{\nu_\theta,n}:=(J_{\nu_\theta}'(j_{\nu_\theta,n}))^2 \left(j_{\nu_\theta,n}\right)^{1+\frac{1}{2k_\theta}}$,
we can now state the main result of this section. 
\begin{Theorem}\label{prop.normder}
Let $\theta\in[1,2)$. Assume $w_0=u_0+iv_0\in L^2(0,1;\mathbb{C})$ and that $\{j_{\nu,n}\}_{n\ge 1}$ is the sequence of positive zeros of $J_\nu$. Then, the following holds: 
\begin{enumerate}
\item[a)]  The solution of \eqref{eq.deg-right} is given by 
\begin{equation}\label{soluz-w}
w^a(x,t)=u^a(x,t)+i v^a(x,t) =e^{(\alpha+i\beta)t}\left(e^{tA}u_0+ie^{tA}v_0\right)(x), 
\end{equation}
where
\begin{equation}\label{eq.sol.u}
\begin{aligned}
 \left(e^{tA}u_0\right) (x) &=\sum_{n = 1}^\infty \dfrac{2e^{-\lambda_n t}}{h_{\nu_\theta,n}} G_{\nu_\theta,n}(x,a) U_n^0(a),
 \end{aligned}
\end{equation}
\begin{equation}\label{eq.sol.v}
\begin{aligned}
 \left(e^{tA}v_0\right) (x) &=\sum_{n = 1}^\infty \dfrac{2e^{-\lambda_n t}}{h_{\nu_\theta,n}} G_{\nu_\theta,n}(x,a) V_n^0(a),
 \end{aligned}
\end{equation}
with
 \begin{equation}\label{eq.U0n}
U^0_n(a) : =  \displaystyle \int_{0}^{j_{\nu_\theta,n}} u_0\left(a+(1-a)\left(\frac{s}{j_{\nu_\theta,n}}\right)^{\frac{1}{k_\theta}}\right) s^{\frac{1}{2k_\theta}} J_{\nu_\theta} (s) \, ds, 
 \end{equation}
 and
  \begin{equation}\label{eq.V0n}
V^0_n(a) : =  \displaystyle \int_{0}^{j_{\nu_\theta,n}} v_0\left(a+(1-a)\left(\frac{s}{j_{\nu_\theta,n}}\right)^{\frac{1}{k_\theta}}\right) s^{\frac{1}{2k_\theta}} J_{\nu_\theta} (s) \, ds\ . 
 \end{equation}
\item[b)] The normal derivative at the boundary is given by
\begin{equation}\label{eq.normderwxa}
  \partial_x w^a(1,t)=  \partial_x u^a(1,t)+i \partial_x v^a(1,t),
\end{equation}
where
 \begin{equation}\label{eq.normderuxa}
\partial_x u^a(1,t)  =  e^{\alpha t}\displaystyle\sum_{n = 1}^\infty \dfrac{2 k_\theta e^{-\lambda_n t}}{(1-a)d_{\nu_\theta,n}}\left[\cos(\beta t)U_n^0(a)-\sin(\beta t)V_n^0(a)\right],
 \end{equation} 
\begin{equation}\label{eq.normdervxa}
\partial_x v^a(1,t)  =  e^{\alpha t}\displaystyle\sum_{n = 1}^\infty \dfrac{2 k_\theta e^{-\lambda_n t}}{(1-a)d_{\nu_\theta,n}}\left[\sin(\beta t)U_n^0(a)+\cos(\beta t)V_n^0(a)\right].
 \end{equation} 
\end{enumerate} 
 Moreover, the map $a\mapsto e^{-\lambda_n t}$ is strictly decreasing for all $t>0$ and $n\ge 1$. 
 \end{Theorem}
   \begin{Remark}\label{regularity}
      Observe that, if $u_0,v_0\in Lip([0,1])$, the normal derivative \eqref{eq.normderwxa} is of class $C^1([\tau,\gamma])$ w.r.t. $a$ for positive time and for every compact interval $[\tau,\gamma]\subset (0,1)$. 
      Indeed, by \eqref{eq.normderuxa} and \eqref{eq.normdervxa}, each general term
    of the series depends $C^1$-smoothly on $a$, since the coefficients
    $
     \begin{pmatrix}
      U_n^0(a)\\
      V_n^0(a)
      \end{pmatrix}
      $
     are $C^1([\tau,\gamma])$ whenever $u_0,v_0\in Lip([0,1])$.
  \end{Remark}
 \noindent 
 \textbf{Proof of Theorem~\ref{prop.normder}:}
\textit{a):} 
The solution $w^a$ to \eqref{eq.deg-right} can be represented as follows: 
 \begin{equation}\label{eq.solpr}
 u^a(x,t)+i v^a(x,t) =e^{(\alpha+i\beta)t}\left(e^{tA}u_0+ie^{tA}v_0\right)(x) = e^{(\alpha+i\beta)t}\displaystyle \sum_{n=1}^\infty e^{-\lambda_n t } \left(u^0_n+iv^0_n\right)  \phi_n(x)  , 
 \end{equation}
where
\begin{equation*}
u^0_n = \displaystyle \int_{a}^1 u_0(x) \phi_n(x) \, dx, \qquad\qquad v^0_n = \displaystyle \int_{a}^1 v_0(x) \phi_n(x) \, dx. 
\end{equation*}
Taking into account \eqref{eq.autofun} and performing the change of variables 
$s = j_{\nu_\theta,n}\left(\frac{x-a}{1-a}\right)^{k_\theta}$, we obtain
 \begin{equation}\label{eq.uon}
\begin{split}
u^0_n &= \displaystyle \int_{a}^1 u_0(x) \phi_n(x) \, dx 
  = 
 \dfrac{\sqrt{2 k_\theta }}{\sqrt{1-a}\, \vert J'_{\nu_\theta} (j_{\nu_\theta,n} ) \vert}   \displaystyle \int_{a}^1 u_0(x) \left(\dfrac{x-a}{1-a}\right)^{\frac{1-\theta}{2}} J_{\nu _\theta} \left(j_{\nu_\theta,n} \left(\dfrac{x-a}{1-a}\right) ^{k_\theta}\right) \, dx  
 \\
 & =  \dfrac{\sqrt{2 k_\theta }}{\sqrt{1-a}\, \vert J'_{\nu_\theta} (j_{\nu_\theta,n} ) \vert}  \displaystyle \int_{0}^{j_{\nu_\theta,n}} u_0\left(a+(1-a)\left(\dfrac{s}{j_{\nu_\theta,n}}\right)^{\frac{1}{k_\theta}} \right)J_{\nu_\theta} (s) \dfrac{1-a}{j_{\nu_\theta,n} k_\theta} \left(\dfrac{s}{j_{\nu_\theta,n}}\right)^{\frac{1}{2k_\theta}} \, ds
  \\
 & : =
     \dfrac{\sqrt{2(1-a)}}{\vert J'_{\nu_\theta} (j_{\nu_\theta,n})\vert \sqrt{k_\theta}\left(j_{\nu_\theta,n}\right)^{1+\frac{1}{2k_\theta}} }  U_n^0(a) 
 \end{split} 
\end{equation}
 and similarly
 \begin{equation}\label{eq.von}
     v^0_n = \displaystyle \int_{a}^1 v_0(x) \phi_n(x) \, dx := \dfrac{\sqrt{2(1-a)}}{\vert J'_{\nu_\theta} (j_{\nu_\theta,n})\vert \sqrt{k_\theta}\left(j_{\nu_\theta,n}\right)^{1+\frac{1}{2k_\theta}} }  V_n^0(a),
 \end{equation}
 where $U_n^0(a)$ and $V_n^0(a)$ are given by \eqref{eq.U0n} and \eqref{eq.V0n}. 
 From \eqref{eq.solpr}, using \eqref{eq.autofun}, \eqref{eq.uon} and \eqref{eq.von}, we deduce  
 \begin{equation*}
     w^a=e^{(\alpha+i\beta)t}\sum_{n = 1}^\infty \dfrac{2e^{-\lambda_n t}}{h_{\nu_\theta,n}} G_{\nu_\theta,n}(x,a) \left(U_n^0+iV_n^0\right).
 \end{equation*}
 \noindent 
 \textit{b)} Taking into account \eqref{soluz-w}, \eqref{eq.sol.u} and \eqref{eq.sol.v}, we get $\partial_x w^a(x,t)=  \partial_x u^a(x,t)+i \partial_x v^a(x,t)$ where
 \begin{equation*}
 \begin{aligned}
  \partial_x u^a(x,t) =e^{\alpha t} \displaystyle\sum_{n = 1}^\infty &\Bigg\{\dfrac{2e^{-\lambda_n t}\left[\cos(\beta t)U_n^0(a)-\sin(\beta t)V_n^0(a)\right]}{h_{\nu_\theta,n}}\cdot\\
  &\quad \cdot\Bigg[ \dfrac{1-\theta}{2(1-a)}\left(\frac{x-a}{1-a} \right)^{-\frac{1+\theta}{2}}J_{\nu_\theta} \left(j_{\nu_\theta,n} \left(\dfrac{x-a}{1-a}\right)^{k_\theta}\right)
  \\
  &\qquad +\left(\frac{x-a}{1-a} \right)^{\frac{1-2\theta}{2}}\dfrac{j_{\nu_\theta,n}k_\theta}{1-a}J'_{\nu_\theta} \left(j_{\nu_\theta,n} \left(\dfrac{x-a}{1-a}\right)^{k_\theta}\right)\Bigg]\Bigg\}.
  \end{aligned}
  \end{equation*}
and
\begin{equation*}
 \begin{aligned}
  \partial_x v^a(x,t) = e^{\alpha t}\displaystyle\sum_{n = 1}^\infty &\Bigg\{\dfrac{2e^{-\lambda_n t}\left[\sin(\beta t)U_n^0(a)+\cos(\beta t)V_n^0(a)\right]}{h_{\nu_\theta,n}}\cdot\\
  &\quad\cdot \Bigg[ \dfrac{1-\theta}{2(1-a)}\left(\frac{x-a}{1-a} \right)^{-\frac{1+\theta}{2}}J_{\nu_\theta} \left(j_{\nu_\theta,n} \left(\dfrac{x-a}{1-a}\right)^{k_\theta}\right)\\
  &\qquad +
  \left(\frac{x-a}{1-a} \right)^{\frac{1-2\theta}{2}}\dfrac{j_{\nu_\theta,n}k_\theta}{1-a}J'_{\nu_\theta} \left(j_{\nu_\theta,n} \left(\dfrac{x-a}{1-a}\right)^{k_\theta}\right)\Bigg]\Bigg\}.
  \end{aligned}
  \end{equation*}
Hence, evaluating for $x = 1$ and exploiting $J_{\nu_\theta} \left(j_{\nu_\theta,n}\right)=0$, we easily obtain \eqref{eq.normderwxa}, \eqref{eq.normderuxa} and \eqref{eq.normdervxa}. 
 
Finally, since 
\begin{equation*}
\partial_a \left(e^{-\lambda_n t}\right) =  -  \dfrac{2k_\theta^3 j_{\nu_\theta,n}^2 t}{(1-a)^{2k_\theta+1}} e^{-\lambda_n t} <0,
\end{equation*}
 we also deduce that the map $a\mapsto e^{-\lambda_n t}$ is strictly decreasing for all $t>0$ and $n\ge 1$.
 \hfill$\blacksquare$

 \section{Lipschitz stability with one point measurement}\label{sec.stablility}
 
Exploiting the explicit expression of the solution given in Theorem \ref{prop.normder}, we will present sufficient conditions for a Lipschitz stability result with a one-point measurement. 
We will also show an example of initial conditions for which a Lipschitz stability estimate can be obtained.

\begin{Theorem}\label{th.stability}
Let $\theta \in [1,2)$ and assume that $u_0,v_0\in Lip([0,1])$. Let $w^{a_1}$ and $w^{a_2}$ be the solutions to \eqref{eq.deg-right},  with $w_0=u_0+iv_0$, corresponding to the degeneracy points $a_1$ and $a_2$, respectively. Assume that there exist $\delta>0$ and $[\tau,\gamma]\subset (0,1)$ such that 
\begin{equation}\label{th2.assump}
 \left|\begin{pmatrix}
    U_1^0(a)\\
    V_1^0(a)
\end{pmatrix}\right|\ge \delta, \qquad \forall\, a\in [\tau,\gamma],
\end{equation}
with $U_1^0(a)$ and $V_1^0(a)$ given by \eqref{eq.U0n} and \eqref{eq.V0n} with $n=1$, respectively.

Then, there exist $t_0(u_0,v_0,\delta,L,\theta,\tau)>0$ and a constant $C>0$ such that the stability estimate 
  \begin{equation}\label{eq.th} 
 |a_2-a_1| \le   C |\partial_x w^{a_2}(1,t) - \partial_x w^{a_1}(1,t)|  
 \end{equation}
 holds 
 \begin{itemize}
     \item for all $a_1,a_2\in [\tau,\gamma]$ and for all $t\in  [t_0, t_1]$ (with $t_1>t_0$), if $\lambda_1(\gamma)>\alpha$;
     \item for all $a_1,a_2\in [\tau,\gamma]$ and for all $t\geq t_0$, if $\lambda_1(\gamma)\leq\alpha$,
 \end{itemize}
 where $\lambda_1(\gamma)=k_\theta^2\dfrac{j_{\nu_\theta,1}^2}{(1-\gamma)^{2k_\theta}}$ and $\alpha$ is the real part of the coefficient $c$ in \eqref{eq.deg}.
\end{Theorem}

 \begin{Remark}\label{rem:stability}
Observe that the assumption \eqref{th2.assump} is satisfied for any $\gamma\in (0,1)$ if $|u_0|>0$ or $|v_0|>0$, respectively, for all $x\in (a,1)$. 
In fact, assuming $|u_0|>0$ for all $x\in (a,1)$, then $\left| U^0_1(a)\right| : =  \left|\displaystyle \int_{0}^{j_{\nu_\theta,1}} u_0\left(a+(1-a)\left(\frac{s}{j_{\nu_\theta,1}}\right)^{\frac{1}{k_\theta}}\right) s^{\frac{1}{2k_\theta}} J_{\nu_\theta} (s) \, ds\,\right| >0$ if $a\in [\tau,\gamma]$, the integrand being strictly positive or strictly negative on $(0,j_{\nu_\theta,1})$.
A similar argument can be made for $\left| V^0_1(a)\right|$. 
 \end{Remark}

 \noindent
 \textbf{Proof of Theorem~\ref{th.stability}:}  
  Using the explicit expression of the normal derivative at the boundary, given by \eqref{eq.normderwxa}, \eqref{eq.normderuxa} and \eqref{eq.normdervxa} and Remark \ref{regularity}, we compute the following:
 \begin{equation}\label{eq1.th2}
 \begin{aligned}
&\left|\partial_a(\partial_x w^a(1,t))\right|  =
\left|\partial_a\begin{pmatrix}
    \partial_x u^a(1,t)\\
    \partial_x v^a(1,t)
\end{pmatrix}\right|=  \left| e^{\alpha t}\displaystyle\sum_{n = 1}^\infty \dfrac{2 k_\theta e^{-\lambda_n t}}{(1-a)d_{\nu_\theta,n}}R(\beta t) 
\begin{pmatrix}
    F^1_n(a)\\
    F^2_n(a)
\end{pmatrix}\right|\\
&=\Bigg| e^{\alpha t}R(\beta t)\Bigg( \dfrac{2 k_\theta e^{-\lambda_1 t}}{(1-a)d_{\nu_\theta,1}} 
\begin{pmatrix}
    F^1_1(a)\\
    F^2_1(a)
\end{pmatrix}+\displaystyle\sum_{n = 2}^\infty \dfrac{2 k_\theta e^{-\lambda_n t}}{(1-a)d_{\nu_\theta,n}}
\begin{pmatrix}
    F^1_n(a)\\
    F^2_n(a)
\end{pmatrix}\Bigg)\Bigg|\ ,
\end{aligned}
\end{equation}
 where 
 \begin{equation}\label{rotation}
      R(\beta t):= \begin{pmatrix}
    \cos(\beta t) & -\sin(\beta t)\\
    \sin(\beta t) & \cos(\beta t)\\
\end{pmatrix},
  \end{equation}
 \begin{equation*}
    F^1_n(a):= \Big[(U_n^0)'(a) +\dfrac{1}{1-a}\left(1 
-2k_\theta^3 \dfrac{j_{\nu_\theta,n}^2 t}{(1-a)^{2k_\theta}}\right) U_n^0(a)\Big],\quad
 \end{equation*}
 \begin{equation*}
    F^2_n(a):= \Big[(V_n^0)'(a)+\dfrac{1}{1-a}\left(1 -2k_\theta^3 \dfrac{j_{\nu_\theta,n}^2 t}{(1-a)^{2k_\theta}}\right) V_n^0(a)\Big] 
 \end{equation*}
 and $U_n^0(a)$, $V_n^0(a)$ are given by \eqref{eq.U0n} and \eqref{eq.V0n} and $\lambda_n$ by \eqref{eigenvalues}.

 Using \eqref{eq1.th2}, we obtain
\begin{equation}
\label{estimate}
\begin{aligned}
    \left|\partial_a(\partial_x w^a(1,t))\right| & 
    \geq e^{\alpha t} e^{-\lambda_1 t}2 k_\theta
    \Bigg(  \dfrac{1}{\vert d_{\nu_\theta,1}\vert} 
\left|\begin{pmatrix}
    F^1_1(a)\\
    F^2_1(a)
\end{pmatrix}\right|-\displaystyle\sum_{n = 2}^\infty \dfrac{ e^{-(\lambda_n-\lambda_1) t}}{\vert d_{\nu_\theta,n}\vert}
\left|\begin{pmatrix}
    F^1_n(a)\\
    F^2_n(a)
\end{pmatrix}\right|\Bigg)\\
& 
    \geq e^{\alpha t} e^{-\lambda_1 t}2 k_\theta
    \Bigg(  \dfrac{1}{\vert d_{\nu_\theta,1}\vert}
    \left|\begin{pmatrix}
    F^1_1(a)\\
    F^2_1(a)
\end{pmatrix}\right|
-\displaystyle\sum_{n = 2}^\infty \dfrac{ e^{-(\lambda_n-\lambda_1) t}\left( \vert F^1_n(a)\vert+\vert F^2_n(a)\vert\right)}{\vert d_{\nu_\theta,n}\vert}
\Bigg).
    \end{aligned}
\end{equation}
Now, let us introduce the notation $M_1:=||u_0||_\infty$, $M_2:=||u'_0||_\infty$, $M_3:=||v_0||_\infty$, $M_4:=||v'_0||_\infty$.

The first term can be estimated in the following way: 
  \begin{equation*}
 \begin{split}
\left|\begin{pmatrix}
    F^1_1(a)\\
    F^2_1(a)
\end{pmatrix}\right|^2&=\left((U_1^0)'(a) +\dfrac{1}{1-a}\left(1-2k_\theta^3 \dfrac{j_{\nu_\theta,1}^2 t}{(1-a)^{2k_\theta}}\right) U_1^0(a)\right)^2\\
&\quad\quad+\left((V_1^0)'(a)+\dfrac{1}{1-a}\left(1-2k_\theta^3 \dfrac{j_{\nu_\theta,1}^2 t}{(1-a)^{2k_\theta}}\right) V_1^0(a)\right)^2
\\
&  = \left((U_1^0)'(a)\right)^2+\left((V_1^0)'(a)\right)^2
   +\dfrac{1}{(1-a)^2}\left(1-2k_\theta^3 \dfrac{j_{\nu_\theta,1}^2 t}{(1-a)^{2k_\theta}}\right)^2 \left( \left(U_1^0(a)\right)^2 + \left(V_1^0(a)\right)^2 \right)
   \\
&\quad +\dfrac{2}{1-a}\left(1-2k_\theta^3 \dfrac{j_{\nu_\theta,1}^2 t}{(1-a)^{2k_\theta}}\right)\left(U_1^0(a)(U_1^0)'(a) +V_1^0(a)(V_1^0)'(a) \right)
\\
& \geq
   \dfrac{1}{2(1-a)^2}\left(1-2k_\theta^3 \dfrac{j_{\nu_\theta,1}^2 t}{(1-a)^{2k_\theta}}\right)^2 \left( \left(U_1^0(a)\right)^2 + \left(V_1^0(a)\right)^2 \right)
   -\left(\left( (U_1^0)'(a)\right)^2 +
 \left( (V_1^0)'(a)\right)^2\right)
   \\
\end{split}
\end{equation*}
\begin{equation}\label{eq2.th2}
    \begin{split}
&  \geq
\dfrac{1}{2}\left(2k_\theta^6 j_{\nu_\theta,1}^4 t^2-1\right)\delta^2
\\
& \quad
-\left(
\displaystyle \int_{0}^{j_{\nu_\theta,1}} u'_0\Bigg(a+(1-a)\left(\frac{s}{j_{\nu_\theta,1}}\right)^{\frac{1}{k_\theta}} \Bigg)\Bigg(1-\left(\dfrac{s}{j_{\nu_\theta,1}}\right)^{\frac{1}{k_\theta}} \Bigg)s^{\frac{1}{2k_\theta}} J_{\nu_\theta} (s) \, ds\right)^2
\\
& \quad
-\left(
\displaystyle \int_{0}^{j_{\nu_\theta,1}} v'_0\Bigg(a+(1-a)\left(\frac{s}{j_{\nu_\theta,1}}\right)^{\frac{1}{k_\theta}} \Bigg)\Bigg(1-\left(\dfrac{s}{j_{\nu_\theta,1}}\right)^{\frac{1}{k_\theta}} \Bigg)s^{\frac{1}{2k_\theta}} J_{\nu_\theta} (s) \, ds\right)^2
\\
& \geq \dfrac{1}{2}\left(2k_\theta^6 j_{\nu_\theta,1}^4 t^2-1\right)\delta^2-\left(M_2^2+M_4^2\right)
\left(\displaystyle \int_{0}^{j_{\nu_\theta,1}} s^{\frac{1}{2k_\theta}} J_{\nu_\theta}(s)\, ds\right)^2
\\
& =\dfrac{1}{2} \left(2k_\theta^6 j_{\nu_\theta,1}^4 t^2-1\right)\delta^2 -j_{\nu_\theta,1}^{2(\nu_\theta+1)} \left(J'_{\nu_\theta}(j_{\nu_\theta,1})\right)^2  \left(M_2^2+M_4^2\right)  \qquad\qquad \forall\,  a\in [\tau,\gamma],
\end{split}
 \end{equation}
where we have exploited \eqref{th2.assump}, the binomial inequality, $\frac{1}{2k_\theta}=\nu_{\theta}+1$ and property d) in Lemma \ref{lemma.Bessel}.
 
 For all $L>0$, we have 
  \begin{equation}\label{eq3.th2}
 \left|\begin{pmatrix}
    F^1_1(a)\\
    F^2_1(a)
\end{pmatrix}\right|^2 \geq L^2\quad \forall\, a\in [\tau,\gamma], \quad
 \forall\, t\ge \bar{t}, 
 \end{equation}
 where  
   \begin{equation*}\label{eqTu0.th2}
\bar{t}(u_0,v_0,\delta,L,\theta) =  \dfrac{1}{\sqrt{2}k_\theta^3 j_{\nu_\theta,1}^2 } \sqrt{1+\dfrac{2}{\delta^2}\left( L^2 + j_{\nu_\theta,1}^{2(\nu_\theta+1)} \left(J'_{\nu_\theta}(j_{\nu_\theta,1})\right)^2  \left(M_2^2+M_4^2\right)\right)} .
 \end{equation*}
\noindent 
Therefore, taking into account \eqref{eq2.th2} and \eqref{eq3.th2}, from \eqref{estimate} we get  
\begin{equation}\label{eq4.th2}
 \begin{split}
\left|\partial_a(\partial_x w^a(1,t))\right| & \ge e^{\alpha t}e^{-\lambda_1 t}2k_\theta \Bigg[ \dfrac{L}{\vert d_{\nu_\theta,1}\vert} 
-
 \displaystyle\sum_{n = 2}^\infty \dfrac{e^{-(\lambda_n -\lambda_1)t}\left( \vert F^1_n(a)\vert+\vert F^2_n(a)\vert\right)}{\vert d_{\nu_\theta,n}\vert}
\Bigg]
\end{split}
 \end{equation}
 for all  $t\ge \bar{t}$.
 So, we have to show that the second term on the right-hand side of \eqref{eq4.th2} is small for $t$ large enough. 
 We have
 \begin{equation*}\label{eq5.th2}
 \begin{split}
\vert F^1_n(a)\vert&=\Big |(U_n^0)'(a) +\dfrac{1}{1-a}\left(1-2k_\theta^3 \dfrac{j_{\nu_\theta,n}^2 t}{(1-a)^{2k_\theta}}\right) U_n^0(a)  \Big| 
\\
&= \Bigg| \displaystyle \int_{0}^{j_{\nu_\theta,n}}\Bigg[  u'_0\Bigg(a+(1-a)\left(\frac{s}{j_{\nu_\theta,n}}\right)^{\frac{1}{k_\theta}} \Bigg)\Bigg(1-\left(\dfrac{s}{j_{\nu_\theta,n}}\right)^{\frac{1}{k_\theta}} \Bigg) \\
& \qquad\qquad\qquad +\dfrac{1}{1-a}\left(1-2k_\theta^3\dfrac{j_{\nu_\theta,n}^2 t}{(1-a)^{2k_\theta}}\right)u_0\Bigg(a+(1-a)\left(\frac{s}{j_{\nu_\theta,n}}\right)^{\frac{1}{k_\theta}} \Bigg)\Bigg]s^{\frac{1}{2k_\theta}} J_{\nu_\theta} (s) \, ds\Bigg| 
\\
& \leq \Big( M_2 +\dfrac{1}{1-a}\left|1- 2k_\theta^3\dfrac{j_{\nu_\theta,n}^2 t}{(1-a)^{2k_\theta}}\right| M_1  \Big)  \int_{0}^{j_{\nu_\theta,n}}s^{\frac{1}{2k_\theta}} |J_{\nu_\theta} (s)| \, ds
 \\
& \leq \dfrac{(j_{\nu_\theta,n})^{\frac{1}{2k_\theta}+1}}{\frac{1}{2k_\theta}+1}K_n^1(t), \qquad \forall\, t\geq \dfrac{(1-a)^{2k_\theta}}{2k_\theta^3 j_{\nu_\theta,n}^2}
\end{split}
 \end{equation*}
where we have taken into account property c) in Lemma \ref{lemma.Bessel} and set
\begin{equation*}
K_n^1(t): = M_2 + \dfrac{1}{1-a}\left(2k_\theta^3\dfrac{j_{\nu_\theta,n}^2 t}{(1-a)^{2k_\theta}}-1\right) M_1 .
\end{equation*}
Similarly, one can obtain the estimate
\begin{equation*}
    \vert F^2_n(a)\vert\leq \dfrac{(j_{\nu_\theta,n})^{\frac{1}{2k_\theta}+1}}{\frac{1}{2k_\theta}+1}K_n^2(t) , \quad \text{where} \quad K_n^2(t): = M_4 +\dfrac{1}{1-a}\left( 2k_\theta^3\dfrac{j_{\nu_\theta,n}^2 t}{(1-a)^{2k_\theta}}-1\right) M_3 .
\end{equation*}
Hence, the second term in the brackets in \eqref{eq4.th2} can be estimated as follows:
 \begin{equation}\label{eq6.th2}
 \begin{aligned}
 &\displaystyle\sum_{n = 2}^\infty \dfrac{e^{-(\lambda_n -\lambda_1)t}\left( \vert F^1_n(a)\vert+\vert F^2_n(a)\vert\right)}{\vert d_{\nu_\theta,n}\vert}
 \le \displaystyle\sum_{n = 2}^\infty \dfrac{e^{-(\lambda_n -\lambda_1)t}}{\vert d_{\nu_\theta,n}\vert}
 \, \dfrac{(j_{\nu_\theta,n})^{\frac{1}{2k_\theta}+1}}{\dfrac{1}{2k_\theta}+1}\left(K_n^1(t)+K_n^2(t)\right)\\
 &= \displaystyle\sum_{n = 2}^\infty \dfrac{e^{-(\lambda_n -\lambda_1)t}}{\vert J_{\nu_\theta}'(j_{\nu_\theta,n})\vert j_{\nu_\theta,n}}
 \, \dfrac{(j_{\nu_\theta,n})^{2}}{\nu_\theta+2}\left(K_n^1(t)+K_n^2(t)\right) : = R_1+R_2.
 \end{aligned}
 \end{equation} 
 Since $|J'_{\nu_\theta} (j_{\nu_\theta,n})|>0$, we deduce that 
 $\lim_{n\to \infty } j_{\nu_\theta,n} |J'_{\nu_\theta} (j_{\nu_\theta,n})| \ge M >0$.
 Therefore, we can estimate the expression of $R_1$ in \eqref{eq6.th2} in the following way:
 \begin{equation}\label{eqR1}
 \begin{split}
 R_1 & \le \dfrac{1}{(\nu_\theta+2)M} \displaystyle\sum_{n = 2}^\infty  j_{\nu_\theta,n}^{2} \, K_n^1(t) \, e^{-(\lambda_n -\lambda_1)t} 
 \\
&  \le
\dfrac{1}{(\nu_\theta+2)M} \Bigg[ \left( M_2-\dfrac{M_1}{1-a}\right) \displaystyle\sum_{n = 2}^\infty  j_{\nu_\theta,n}^{2} e^{-(\lambda_n -\lambda_1)t}   
+
 2k_\theta^3\dfrac{t  M_1}{(1-a)^{2k_\theta+1}} \displaystyle\sum_{n = 2}^\infty j_{\nu_\theta,n}^{4} e^{-(\lambda_n -\lambda_1)t}   \Bigg].
 \end{split}
 \end{equation} 
In addition, 
\begin{equation*}
j_{\nu_\theta,n}^{2} e^{-\lambda_n t} \le e^{j_{\nu_\theta,n}^{2}} e^{-\frac{k_\theta^2 t }{(1-a)^{2k_\theta}}j_{\nu_\theta,n}^{2}} \le e^{\frac{(1-a)^{2k_\theta}-k_\theta^2t}{(1-a)^{2k_\theta}}j_{\nu_\theta,n}^{2}}
\le 
e^{-\frac{k_\theta^2t }{2(1-a)^{2k_\theta}}j_{\nu_\theta,n}^{2}},\qquad \forall\, t\ge \frac{2(1-a)^{2k_\theta}}{k_\theta^2}
\end{equation*} 
and, using $x^2\le e^x$ $\forall\, x\ge 0$, we get 
\begin{equation*}
j_{\nu_\theta,n}^{4} e^{-\lambda_n t} \le e^{j_{\nu_\theta,n}^{2}} e^{-\frac{k_\theta^2 t }{(1-a)^{2k_\theta}}j_{\nu_\theta,n}^{2}} \le 
e^{-\frac{k_\theta^2 t }{2(1-a)^{2k_\theta}}j_{\nu_\theta,n}^{2}},\qquad \forall\, t\ge \frac{2(1-a)^{2k_\theta}}{k_\theta^2}.
\end{equation*}
\noindent
By considering these estimates, we conclude that
 \begin{equation*}\label{eq8.th2}
 \begin{split}
 R_1 &  \le
\dfrac{1}{(\nu_\theta+2)M}  \Bigg[   M_2-\dfrac{M_1}{1-a}    
+
 \dfrac{2k_\theta^3  M_1 t}{(1-a)^{2k_\theta+1}}    \Bigg]
 \displaystyle\sum_{n = 2}^\infty   e^{\left(\lambda_1-\frac{k_\theta^2 j_{\nu_\theta,n}^{2} }{2(1-a)^{2k_\theta}}\right)t}, \qquad  \forall\,  t\ge \frac{2(1-a)^{2k_\theta}}{k_\theta^2}.
 \end{split}
 \end{equation*}
\noindent 
We now claim that
 \begin{equation}\label{eqR1-lim}
\lim_{t\to +\infty} \dfrac{1}{(\nu_\theta+2)M}  \Bigg[   M_2-\dfrac{M_1}{1-a}    
+
 \dfrac{2k_\theta^3  M_1 t}{(1-a)^{2k_\theta+1}}    \Bigg]
 \displaystyle\sum_{n = 2}^\infty   e^{\left(\lambda_1-\frac{k_\theta^2 j_{\nu_\theta,n}^{2} }{2(1-a)^{2k_\theta}}\right)t}  = 0. 
 \end{equation}
Regarding the last term in \eqref{eqR1-lim}, for $n$ sufficiently large so that 
    $\lambda_1-\frac{k_\theta^2 j_{\nu_\theta,n}^{2} }{2(1-a)^{2k_\theta}} \leq
    -\frac{k_\theta^3 j_{\nu_\theta,n}^{2} }{2(1-a)^{2k_\theta}}$,
we have that
  \begin{equation}\label{eq11.th2}
t e^{\left(\lambda_1-\frac{k_\theta^2 j_{\nu_\theta,n}^{2} }{2(1-a)^{2k_\theta}}\right)t} 
 \le 
 t e^{-\frac{k_\theta^3 j_{\nu_\theta,n}^{2}t }{2(1-a)^{2k_\theta}}} \frac{k_\theta^3 j_{\nu_\theta,n}^{2} }{2(1-a)^{2k_\theta}} \frac{2(1-a)^{2k_\theta} }{k_\theta^3 j_{\nu_\theta,n}^{2}} \le \frac{2(1-a)^{2k_\theta} }{ek_\theta^3 j_{\nu_\theta,n}^{2}},
 \end{equation} 
 which is summable for $j_{\nu_\theta,n}^2\sim \pi^2n^2$ (see properties f) and g) of Lemma \ref{lemma.Bessel}). 
 Let us remark that in \eqref{eq11.th2} we have used the fact that $  t e^{-\frac{k_\theta^3 j_{\nu_\theta,n}^{2}t }{2(1-a)^{2k_\theta}}} \frac{k_\theta^3 j_{\nu_\theta,n}^{2} }{2(1-a)^{2k_\theta}} \le e^{-1}$, which follows from the standard inequality 
 $xe^{-x} \le e^{-1}$.
  Since $\lim_{t\to \infty}t  e^{\left(\lambda_1-\frac{k_\theta^2 j_{\nu_\theta,n}^{2} }{2(1-a)^{2k_\theta}}\right)t} =0$, Weierstrass M-test (or Lebesgue's dominated convergence theorem applied to the counting measure) yields
$$\lim_{t\to \infty} \displaystyle\sum_{n = 2}^\infty   t  e^{\left(\lambda_1-\frac{k_\theta^2 j_{\nu_\theta,n}^{2} }{2(1-a)^{2k_\theta}}\right)t}=0.$$
With an analogous procedure, we also have  the following
\begin{equation*}
    \lim_{t\to \infty} \displaystyle\sum_{n = 2}^\infty    e^{\left(\lambda_1-\frac{k_\theta^2 j_{\nu_\theta,n}^{2} }{2(1-a)^{2k_\theta}}\right)t}=0.
\end{equation*}

 Similarly, we can obtain an analogous estimate for $R_2$ and claim that
 \begin{equation}\label{eqR2}
 \begin{split}
 R_2 \le
\dfrac{1}{(\nu_\theta+2)M} \Bigg[   \left(M_4-\dfrac{M_3}{1-a}\right) \displaystyle\sum_{n = 2}^\infty  j_{\nu_\theta,n}^{2} e^{-(\lambda_n -\lambda_1)t}   
+
 2k_\theta^3\dfrac{t  M_3}{(1-a)^{2k_\theta+1}} \displaystyle\sum_{n = 2}^\infty j_{\nu_\theta,n}^{4} e^{-(\lambda_n -\lambda_1)t}   \Bigg]
 \end{split}
 \end{equation}
 and
\begin{equation}\label{eqR2-lim}
\lim_{t\to +\infty} \dfrac{1}{(\nu_\theta+2)M}  \Bigg[   M_4-\dfrac{M_3}{1-a}    
+
 \dfrac{2k_\theta^3  M_3 t}{(1-a)^{2k_\theta+1}}    \Bigg]
 \displaystyle\sum_{n = 2}^\infty   e^{\left(\lambda_1-\frac{k_\theta^2 j_{\nu_\theta,n}^{2} }{2(1-a)^{2k_\theta}}\right)t}  = 0. 
 \end{equation}
 Taking into account \eqref{eqR1-lim} and \eqref{eqR2-lim}, we deduce from \eqref{eq4.th2}, \eqref{eqR1} and \eqref{eqR2} that there exists $t_0(u_0,v_0,\delta,L,\theta,\tau)>0$ such that
\begin{equation*}\label{eq12.th2}
|\partial_a(\partial_x w^a(1,t))| \ge e^{\alpha t}e^{-\lambda_1(\gamma) t}2k_\theta  \dfrac{L}{\vert d_{\nu_\theta,1}\vert},   \quad \forall\, t\ge t_0(u_0,v_0,\delta,L,\theta,\tau).
\end{equation*}
To conclude, let us obtain the stability estimate. 
For all $a_1,a_2\in [\tau,\gamma]$, we get
 \begin{equation*}
| \partial_x w^{a_2}(1,t) - \partial_x w^{a_1}(1,t)| \ge e^{\alpha t}e^{-\lambda_1(\gamma) t}2k_\theta  \dfrac{L}{\vert d_{\nu_\theta,1}\vert} |a_2-a_1|,\quad \forall\, t\ge t_0(u_0,v_0,\delta,L,\theta,\tau)\ .
 \end{equation*}
 If $\lambda_1(\gamma)>\alpha$, by fixing $t_1>t_0$, we obtain \eqref{eq.th} with $C = e^{(\lambda_1(\gamma)-\alpha) t_1}  \dfrac{\vert d_{\nu_\theta,1}\vert}{2k_\theta L} $.
 If $\lambda_1(\gamma)\leq\alpha$, we get \eqref{eq.th} with $C = e^{(\lambda_1(\gamma)-\alpha) t_0}  \dfrac{\vert d_{\nu_\theta,1}\vert}{2k_\theta L} $.
 This ends the proof. \hfill $\blacksquare$
\vskip 2pt

In Theorem \ref{th.stability}, we assume $a\in[\tau,\gamma]$, which is a compact interval that excludes points 0 and~1. 
The exclusion of the right endpoint is due to the specific point where we perform the normal derivative measurements. 
The exclusion of zero, on the other hand, is simply because we are considering a Dirichlet boundary condition at zero.

However, if we remove this latter condition and analyze the solutions of the problem within the interval $(a,1)$, we can consider a compact interval of the form $[0,\gamma]$, allowing for a degeneracy at the left boundary.
By repeating the proof of Theorem \ref{th.stability}, we can prove the following result:
\begin{Theorem}
Let $\theta \in [1,2)$ and assume that $u_0,v_0\in Lip([0,1])$. Let $w^{a_1}$ and $w^{a_2}$ be the solutions to \eqref{eq.deg-right},  with $w_0=u_0+iv_0$, corresponding to the degeneracy points $a_1$ and $a_2$, respectively, with $0\leq a_1,a_2\leq\gamma<1$. Assume that there exist $\delta>0$ and $[0,\gamma]\subset [0,1)$ such that 
\begin{equation*}
 \left|\begin{pmatrix}
    U_1^0(a)\\
    V_1^0(a)
\end{pmatrix}\right|\ge \delta, \qquad \forall\, a\in [0,\gamma],
\end{equation*}
with $U_1^0(a)$ and $V_1^0(a)$ given by \eqref{eq.U0n} and \eqref{eq.V0n} with $n=1$, respectively.

Then, there exist $t_0(u_0,v_0,\delta,L,\theta,\tau)>0$ and a constant $C>0$ such that the stability estimate 
  \begin{equation*}
 |a_2-a_1| \le   C |\partial_x w^{a_2}(1,t) - \partial_x w^{a_1}(1,t)|  
 \end{equation*}
 holds 
 \begin{itemize}
     \item for all $a_1,a_2\in [0,\gamma]$ and for all $t\in  [t_0, t_1]$ (with $t_1>t_0$), if $\lambda_1(\gamma)>\alpha$;
     \item for all $a_1,a_2\in [0,\gamma]$ and for all $t\geq t_0$, if $\lambda_1(\gamma)\leq\alpha$,
 \end{itemize}
 where $\lambda_1(\gamma)=k_\theta^2\dfrac{j_{\nu_\theta,1}^2}{(1-\gamma)^{2k_\theta}}$ and $\alpha$ is the real part of the coefficient $c$ in \eqref{eq.deg}.
\end{Theorem}

\subsection{Example of admissible initial data for stability estimates}

Now, we will analyze an example of admissible initial data for stability estimates, using a one-point measurement.

We consider the system \eqref{eq.deg-right} with $\theta=1.3$, $\alpha=1$ and $\beta=1/2$. 
We also assume that $u_0(x)=0$ and $v_0(x)=1$, which implies the validity of hypothesis \eqref{th2.assump} (see Remark \ref{rem:stability}).
Following the proof of Theorem \ref{th.stability}, we want to verify if 
\begin{equation}\label{cond-der}
    |\partial_a(\partial_x w^a(1,t))|>0
\end{equation} for $t$ large enough and $\forall a\in[\tau,\gamma]\subset (0,1)$.

\begin{figure}[h!]
\begin{minipage}[t]{0.49\linewidth}
\noindent
 \centering 
 \includegraphics[width=0.98\linewidth]{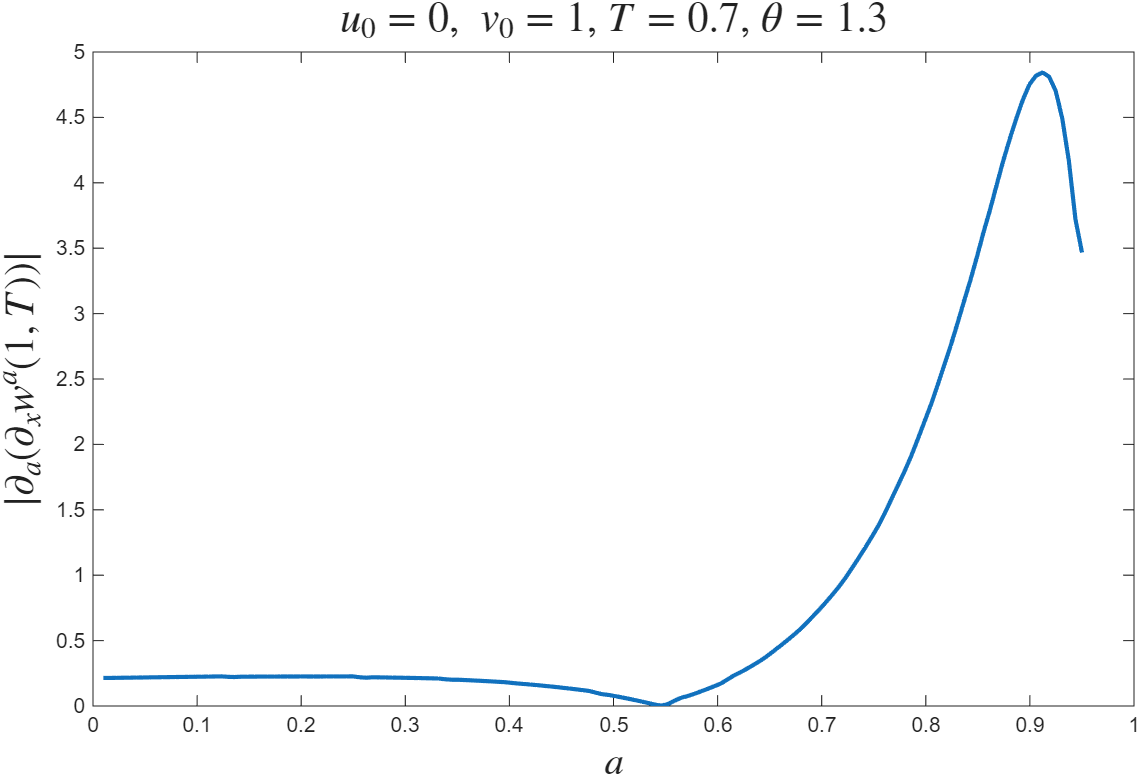}
 \caption{Lack of stability, $T = 0.7$.}
 \label{nostab}
\end{minipage}
\hfill
\begin{minipage}[t]{0.49\linewidth}
\noindent
 \centering 
 \includegraphics[width=0.98\linewidth]{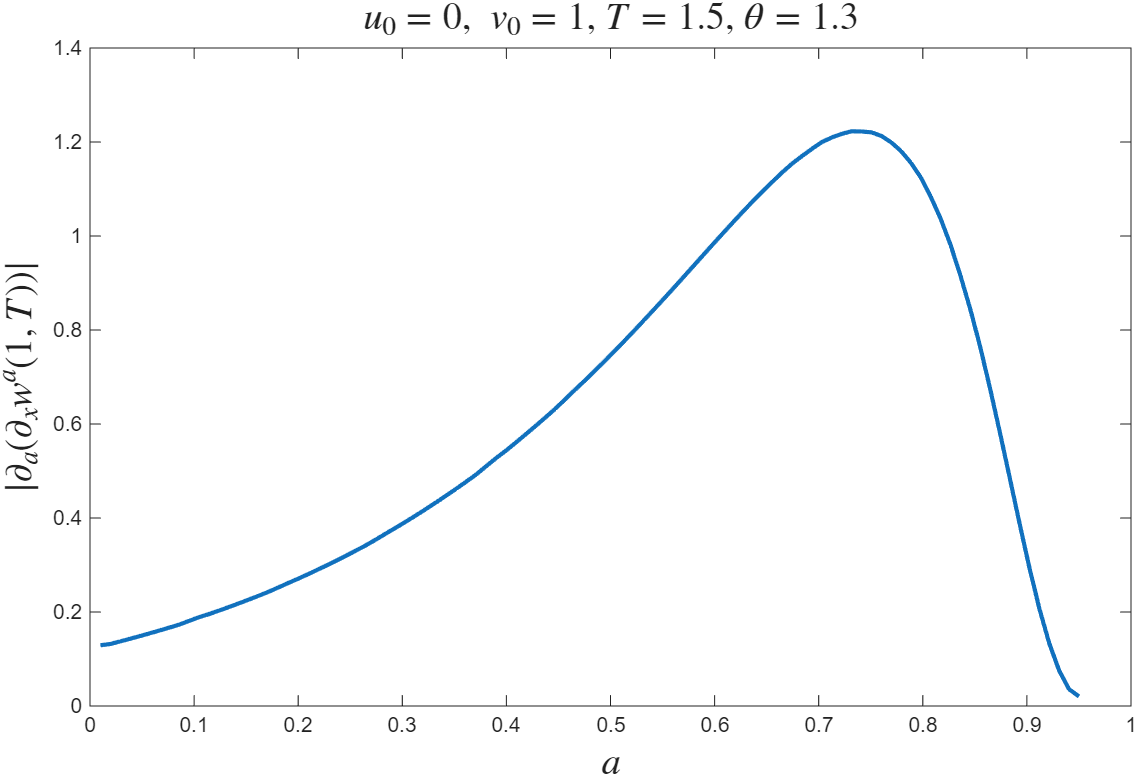}
 \caption{Stability for $t$ large, $T = 1.5$.}
 \label{stab}
\end{minipage}
 \end{figure}

 \begin{figure}
    \centering
    \includegraphics[width=0.5\linewidth]{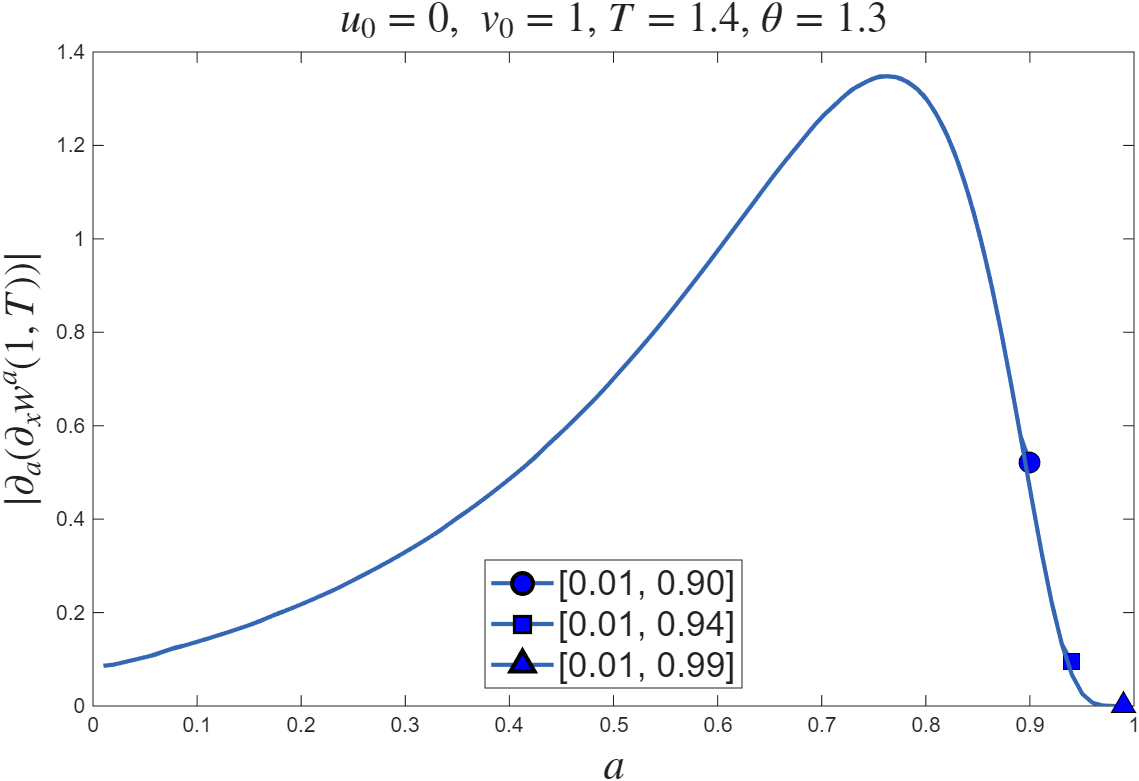}
    \caption{Comparison between several compact intervals $[\tau,\gamma]$.}
    \label{intervals}
\end{figure}

In Figure \ref{nostab}, by fixing a time $T=0.7$, we can see that there exists a point $a$ for which condition \eqref{cond-der} is violated and we cannot guarantee a Lipschitz stability estimate.
Instead, fixing a time large enough (for instance $T=1.5$), we obtain the validity of the condition \eqref{cond-der} $\forall a\in[\tau,\gamma]\subset (0,1)$ (see Figure \ref{stab}). 
We observe that we can get the validity of the condition \eqref{cond-der} as long as $\gamma$ is not too close to the point $1$.
This is due to the fact that it is possible to obtain the result by staying away from the measurement point $x=1$, as we can see in Figure \ref{intervals}, where  we make a comparison considering several compact intervals with different endpoints $\gamma$ (by fixing $T=1.4$ only to make the graphic representation clearer).

\section{Uniqueness results for \textquotedblleft distributed\textquotedblright measurements}\label{sec:uniqueness}

In this section, we present sufficient conditions for the uniqueness result of three more general inverse problems. 
In the first one, we consider the degeneracy point, initial data, and coefficient $c$ as unknowns; in the second one, only the degeneracy point and the initial data are unknown.
Moreover, the third inverse problem concerns the case where the unknowns are represented by the degeneracy power and the initial data.
Unlike in the previous section, where we considered point-wise measurements of $\partial_x w(1,t)$, here we require measurements distributed over a time interval.
For the first inverse problem, in addition to distributed measurements of $\partial_x w(1,t)$, the uniqueness of the coefficient $c$ also requires distributed measurements of $\partial_x w(0,t)$, which also allow us to achieve the uniqueness of the initial data over the entire interval $(0,1)$.
The second and third inverse problems simply require distributed measurements of $\partial_x w(1,t)$, but the uniqueness of the initial data from the $\partial_x w(1,t)$ measurements is confined to the right sub-interval $(a,1)$.

\medskip

We now consider the two sub-problems
\begin{equation}\label{eq.deg-right-uniq}
\left\{\begin{array}{ll}
\partial_t w- \partial_x((x-a)^\theta \partial_x w) -cw= 0, & (x,t)\in(a,1)\times(0,T),  \\[1mm]
 (x-a)^\theta \partial_x w|_{x=a}=0, \quad w(1,t) = 0, &  t\in (0,T),  \\[1mm]
w(x,0) = w_0(x),  & x\in(a,1)
\end{array}\right.
\end{equation}
and
\begin{equation}\label{eq.deg-left-uniq}
\left\{\begin{array}{ll}
\partial_t w- \partial_x((a-x)^\theta \partial_x w) -cw= 0, & (x,t)\in(0,a)\times(0,T),  \\[1mm]
 w(0,t) = 0, \quad (a-x)^\theta \partial_x w|_{x=a}=0, &  t\in (0,T),  \\[1mm]
w(x,0) = w_0(x),  & x\in(0,a),
\end{array}\right.
\end{equation}
where we take into account different unknowns, depending on the inverse problem we are analyzing.
The solution and the normal derivative for the left sub-problem \eqref{eq.deg-left-uniq} can be computed with an analysis similar to that of Section \ref{sec:normal-der}, obtaining
\begin{equation*}
    w(x,t)=e^{(\alpha+i\beta)t}\sum_{n = 1}^\infty \dfrac{2e^{-\lambda_n t}}{h_{\nu_\theta,n}} \left(\frac{a-x}{a} \right)^{\frac{1-\theta}{2}} J_{\nu_\theta} \left(j_{\nu_\theta,n} \left(\dfrac{a-x}{a}\right)^{k_\theta}\right) \left(U_n^0(a)+iV_n^0(a)\right),
\end{equation*}
where 
 \begin{equation}\label{eq.U0n-left}
U^0_n(a) : =  \displaystyle \int_{0}^{j_{\nu_\theta,n}} u_0\left(a-a\left(\frac{s}{j_{\nu_\theta,n}}\right)^{\frac{1}{k_\theta}}\right) s^{\frac{1}{2k_\theta}} J_{\nu_\theta} (s) \, ds, 
 \end{equation}
 and
  \begin{equation}\label{eq.V0n-left}
V^0_n(a) : =  \displaystyle \int_{0}^{j_{\nu_\theta,n}} v_0\left(a-a\left(\frac{s}{j_{\nu_\theta,n}}\right)^{\frac{1}{k_\theta}}\right) s^{\frac{1}{2k_\theta}} J_{\nu_\theta} (s) \, ds 
 \end{equation}
 and
\begin{equation*}
\lambda_n=k_\theta^2\dfrac{j_{\nu_\theta,n}^2}{a^{2k_\theta}}.
\end{equation*}
Moreover, the vector of normal derivatives at the boundary becomes
\begin{equation}\label{normalder-left}
    \begin{pmatrix}
        \partial_x u(0,t)\\
        \partial_x v(0,t)
    \end{pmatrix}=e^{\alpha t}\displaystyle\sum_{n = 1}^\infty \dfrac{-2 k_\theta e^{-\lambda_n t}}{a\, d_{\nu_\theta,n}} R(\beta t)
\begin{pmatrix}
    U_n^0(a)\\
    V_n^0(a)
\end{pmatrix}.
\end{equation}
We can now state the uniqueness result.
 \begin{Theorem}\label{theorem2}
Let $\theta \in [1,2)$, $0<a_1,a_2 <1$ and $0<t_1<t_2$. 
Let $w^{a_1}$ and $w^{a_2}$ be two solutions to \eqref{eq.deg-right-uniq} or \eqref{eq.deg-left-uniq}, corresponding to the initial values $w_0=u_0+iv_0$ and $\widetilde{w}_0=\widetilde{u}_0+i\widetilde{v}_0$, respectively, and the coefficients $c=\alpha+i\beta$ and $\tilde{c}=\tilde{\alpha}+i\tilde{\beta}$, respectively. Assume that 
\begin{equation}\label{th2.eq1}
\begin{pmatrix}
    U_1^0(a_1) \\
    V_1^0(a_1)
\end{pmatrix}\neq
  \begin{pmatrix}
    0 \\
    0
\end{pmatrix}
 \qquad \text{and} \qquad 
 \begin{pmatrix}
    \widetilde{U}_1^0(a_2) \\
    \widetilde{V}_1^0(a_2)
\end{pmatrix}\neq
  \begin{pmatrix}
    0 \\
    0
\end{pmatrix}, 
\end{equation}
where $U^0_1(a_1)$ and $\widetilde{U}^0_1(a_2)$, for $u_0$ and $\widetilde{u}_0$, are given by~\eqref{eq.U0n} or \eqref{eq.U0n-left} (with $n=1$) for \eqref{eq.deg-right-uniq} and \eqref{eq.deg-left-uniq}, respectively; $V^0_1(a_1)$ and $\widetilde{V}^0_1(a_2)$, for $v_0$ and $\widetilde{v}_0$, are given by~\eqref{eq.V0n} or \eqref{eq.V0n-left} (with $n=1$) for \eqref{eq.deg-right-uniq} and \eqref{eq.deg-left-uniq}, respectively.
 
Then $\partial_x w^{a_1}(1,t) = \partial_x w^{a_2}(1,t)$ and $\partial_x w^{a_1}(0,t) = \partial_x w^{a_2}(0,t)$ for $t_1<t<t_2$ imply that $\alpha=\tilde{\alpha}$, $a_1=a_2$, $u_0 = \widetilde{u}_0$, $v_0 = \widetilde{v}_0$ and $\beta=\tilde{\beta}$. 
 \end{Theorem}
\begin{Remark}
    Let us observe that an analogous result can also be obtained assuming that the first~$m-1$ vectors $\begin{pmatrix}
    U_n^0(a_1) \\
    V_n^0(a_1)
\end{pmatrix}$
 and 
 $\begin{pmatrix}
    \widetilde{U}_n^0(a_2) \\
    \widetilde{V}_n^0(a_2)
\end{pmatrix} $ are identically zero and replacing hypothesis \eqref{th2.eq1} with the condition
$\begin{pmatrix}
    U_m^0(a_1) \\
    V_m^0(a_1)
\end{pmatrix}\neq
  \begin{pmatrix}
    0 \\
    0
\end{pmatrix}$
  and 
 $\begin{pmatrix}
    \widetilde{U}_m^0(a_2) \\
    \widetilde{V}_m^0(a_2)
\end{pmatrix}\neq
  \begin{pmatrix}
    0 \\
    0
\end{pmatrix}$.
\end{Remark}

\textbf{Proof of Theorem~\ref{theorem2}:} Note that the functions $t\mapsto \partial_x u^a(1,t)$, $t\mapsto \partial_x v^a(1,t)$, $t\mapsto \partial_x u^a(0,t)$ and $t\mapsto \partial_x v^a(0,t)$ are analytic for all $t>0$. 
Let us consider the left sub-problem \eqref{eq.deg-left-uniq} and set 
\begin{equation*}\label{lambdamu}
\lambda_n:= k_\theta ^2 \frac{j_{\nu _\theta ,n}^2}{a_1^{2k_\theta}}, \quad \mu_n:= k_\theta ^2 \frac{j_{\nu _\theta ,n}^2}{a_2^{2k_\theta}}, \quad n\in \mathbb{N},
\end{equation*}
for which we have $\lambda_1<\lambda_2<\ldots$ and $\mu_1<\mu_2<\ldots$. 

Assume that $a_1\neq a_2$. Without loss of generality, we can assume that $a_1<a_2$, then $\lambda_n>\mu_n$ for $n\in \mathbb{N}$.
Thus, taking into account \eqref{eq.U0n-left}, \eqref{eq.V0n-left} and the condition $\partial_x w^{a_1}(0,t) = \partial_x w^{a_2}(0,t)$ and setting
\begin{equation*}
     W^0_n(a_1):=\begin{pmatrix} 
    U_n^0(a_1)\\
    V_n^0(a_1)
\end{pmatrix}
\quad \text{and} \quad \widetilde{W}^0_n(a_2):=\begin{pmatrix}
    \widetilde{U}^0_n(a_2)\\
    \widetilde{V}^0_n(a_2)
\end{pmatrix},
\end{equation*}
we get
\begin{equation*}
    e^{\alpha t}\displaystyle\sum_{n = 1}^\infty \dfrac{ e^{-\lambda_n t}}{a_1 d_{\nu_\theta,n}} R(\beta t) W^0_n(a_1)
=e^{\tilde{\alpha} t}\displaystyle\sum_{n = 1}^\infty \dfrac{ e^{-\mu_n t}}{a_2 d_{\nu_\theta,n}} R(\tilde{\beta} t) \widetilde{W}^0_n(a_2),
\end{equation*}
which implies
\begin{equation*}
\Bigg| \dfrac{e^{-(\lambda_1-\mu_1) t}}{a_1 d_{\nu_\theta,1}}
W^0_1(a_1)+\displaystyle\sum_{n = 2}^\infty \dfrac{e^{-(\lambda_n-\mu_1) t}}{a_1 d_{\nu_\theta,n}}
W^0_n(a_1)\Bigg|
=\Bigg|e^{(\tilde{\alpha}-\alpha)t}\left(\dfrac{1}{a_2 d_{\nu_\theta,1}} 
\widetilde{W}^0_1(a_2)+\displaystyle\sum_{n = 2}^\infty \dfrac{e^{-(\mu_n-\mu_1) t}}{a_2 d_{\nu_\theta,n}} 
\widetilde{W}^0_n(a_2)\right)\Bigg|.
\end{equation*}
Since $\lambda_n>\mu_n$, letting $t\to\infty$ implies
 $   \Bigg|\dfrac{e^{(\tilde{\alpha}-\alpha)t}}{a_2 d_{\nu_\theta,1}} 
\widetilde{W}^0_1(a_2)\Bigg|\to 0$,
which is true only if $\tilde{\alpha}<\alpha$, due to hypothesis \eqref{th2.eq1}.

Now, let us consider the right sub-problem \eqref{eq.deg-right-uniq} and set 
\begin{equation*}\label{lambdamu}
\lambda_n:= k_\theta ^2 \frac{j_{\nu _\theta ,n}^2}{(1-a_1)^{2k_\theta}}, \quad \mu_n:= k_\theta ^2 \frac{j_{\nu _\theta ,n}^2}{(1-a_2)^{2k_\theta}}, \quad n\in \mathbb{N},
\end{equation*}
for which we have $\lambda_1<\lambda_2<\ldots$ and $\mu_1<\mu_2<\ldots$.  
Due to condition $a_1<a_2$, in this case we have $\lambda_n<\mu_n$ for $n\in \mathbb{N}$.
Hence, taking into account \eqref{eq.U0n}, \eqref{eq.V0n} and the condition $\partial_x w^{a_1}(1,t) = \partial_x w^{a_2}(1,t)$, we get
\begin{equation*}
    e^{\alpha t}\displaystyle\sum_{n = 1}^\infty \dfrac{ e^{-\lambda_n t}}{(1-a_1)d_{\nu_\theta,n}} R(\beta t)
W^0_n(a_1)
=e^{\tilde{\alpha} t}\displaystyle\sum_{n = 1}^\infty \dfrac{ e^{-\mu_n t}}{(1-a_2)d_{\nu_\theta,n}} R(\tilde{\beta} t)
\widetilde{W}^0_n(a_2),
\end{equation*}
which implies
\begin{equation*}
    \Bigg|e^{(\alpha-\tilde{\alpha})t}\left( \dfrac{W^0_1(a_1)}{(1-a_1)d_{\nu_\theta,1}}
+\displaystyle\sum_{n = 2}^\infty \dfrac{e^{-(\lambda_n-\lambda_1) t}}{(1-a_1)d_{\nu_\theta,n}}
W^0_n(a_1)\right)\Bigg|
=\Bigg|\dfrac{e^{-(\mu_1-\lambda_1) t}}{(1-a_2)d_{\nu_\theta,1}} 
\widetilde{W}^0_1(a_2)+\displaystyle\sum_{n = 2}^\infty \dfrac{e^{-(\mu_n-\lambda_1) t}}{(1-a_2)d_{\nu_\theta,n}} 
\widetilde{W}^0_n(a_2)\Bigg|.
\end{equation*}
Since $\lambda_n<\mu_n$, letting $t\to\infty$ implies
$
    \Bigg|\dfrac{e^{(\alpha-\tilde{\alpha})t}}{(1-a_1)d_{\nu_\theta,1}} 
W^0_1(a_1)\Bigg|\to 0$,
which, due to hypothesis \eqref{th2.eq1}, is true only if $\alpha<\tilde{\alpha}$, in contradiction with the previous conclusion. Thus, we can deduce $\alpha=\tilde{\alpha}$.

Using, for example, the condition $\partial_x w^{a_1}(1,t) = \partial_x w^{a_2}(1,t)$, we now prove that $a_1$=$a_2$.
 Therefore, we obtain 
\begin{equation}\label{eq2.uniq2t2}
\displaystyle\sum_{n = 1}^\infty \dfrac{ e^{-\lambda_n t}}{(1-a_1)d_{\nu_\theta,n}} R(\beta t)
W^0_n(a_1)
=\displaystyle\sum_{n = 1}^\infty \dfrac{ e^{-\mu_n t}}{(1-a_2)d_{\nu_\theta,n}} R(\tilde{\beta} t)
\widetilde{W}^0_n(a_2)
\end{equation}
for $t>t_1$.
 Therefore, taking the absolute value, we have 
\begin{equation*}
   \left| \dfrac{1}{(1-a_1)d_{\nu_\theta,1}}
W^0_1(a_1)+\displaystyle\sum_{n = 2}^\infty \dfrac{e^{-(\lambda_n-\lambda_1) t}}{(1-a_1)d_{\nu_\theta,n}}
W^0_n(a_1)\right|
=\left|\dfrac{e^{-(\mu_1-\lambda_1) t}}{(1-a_2)d_{\nu_\theta,1}} 
\widetilde{W}^0_1(a_2)+\displaystyle\sum_{n = 2}^\infty \dfrac{e^{-(\mu_n-\lambda_1) t}}{(1-a_2)d_{\nu_\theta,n}}
\widetilde{W}^0_n(a_2)\right|.
\end{equation*}
Since $\mu_n>\lambda_n$ for all $n\geq1$, we let $t\to\infty$ in the previous equality to get
$
\left|
    W^0_1(a_1)\right|=
    0 $,
in contrast with \eqref{th2.eq1}. Thus, $a_1=a_2$ follows, and this implies $\lambda_n=\mu_n$ for all $n\in \mathbb{N}$ and both problems \eqref{eq.deg-right-uniq} and \eqref{eq.deg-left-uniq}.

Let us now see that $u_0 = \widetilde{u}_0$ and $v_0 = \widetilde{v}_0$. 
Since $\lambda_n=\mu_n$ for all $n\in \mathbb{N}$ and $a_1=a_2\equiv a$, from the absolute value of \eqref{eq2.uniq2t2}, we obtain 
\begin{equation}\label{eq.u0}
\left|\displaystyle\sum_{n = 1}^\infty \dfrac{ e^{-\lambda_n t}}{(1-a)d_{\nu_\theta,n}} 
\left(W^0_n(a)-\widetilde{W}^0_n(a)\right)\right|= 
0
, \quad t>t_1. 
\end{equation}
Let us set $
n_0 = \inf\{ n\ge 1:  W^0_n(a)\neq
    \widetilde{W}^0_n(a)  \}$.
We are going to show that this is an empty set or, equivalently, $n_0=\infty$. Suppose $n_0<\infty$ and multiply the equality \eqref{eq.u0} by $e^{\lambda_{n_0}t}$ to obtain 
\begin{equation*}
\left|\dfrac{ 1}{(1-a)d_{\nu_\theta,n_0}} 
\left(W^0_{n_0}(a)-\widetilde{W}^0_{n_0}(a)\right)+\displaystyle\sum_{n = n_0+1}^\infty \dfrac{ e^{-(\lambda_n-\lambda_{n_0}) t}}{(1-a)d_{\nu_\theta,n}} 
\left(W^0_n(a)-\widetilde{W}^0_n(a)\right)\right|= 
0, 
\end{equation*}
for $t>t_1$.
\noindent
We let $t\to +\infty$ and deduce from the previous equality that 
$W^0_{n_0}(a)=\widetilde{W}^0_{n_0}(a)$,
in contrast with the definition of $n_0$. Therefore, $n_0=\infty$ and 
$W^0_{n}(a)=\widetilde{W}^0_{n}(a)$ $\forall\, n\geq 1$.
From \eqref{soluz-w}, \eqref{eq.sol.u} and \eqref{eq.sol.v}, we conclude that $u_0 =\widetilde{u}_0$ and $v_0 =\widetilde{v}_0$ for $x\in(a,1)$, by the coincidence of all Fourier coefficients. A similar analysis can be performed for the left sub-problem, in order to obtain $u_0 =\widetilde{u}_0$ and $v_0 =\widetilde{v}_0$ for $x\in(0,a)$.

Let us now see that $\beta=\tilde{\beta}$.
From \eqref{eq2.uniq2t2}, we obtain
\begin{equation*}
  R(\beta t)\left(   \dfrac{1}{d_{\nu_\theta,1}}
W^0_1(a)+\displaystyle\sum_{n = 2}^\infty \dfrac{e^{-(\lambda_n-\lambda_1) t}}{d_{\nu_\theta,n}}
W^0_n(a)\right)=R(\tilde{\beta} t)\left(\dfrac{1}{d_{\nu_\theta,1}} 
W^0_1(a)+\displaystyle\sum_{n = 2}^\infty \dfrac{e^{-(\lambda_n-\lambda_1) t}}{d_{\nu_\theta,n}} 
W^0_n(a)\right).
\end{equation*}
Since $\lambda_n>\lambda_1$ for $n\geq 2$, letting $t\to\infty$ implies $\beta=\tilde{\beta}$.
This ends the proof. 
\hfill$\blacksquare$

\begin{Remark} 
Note that assumption \eqref{th2.eq1} is satisfied if
$|w_0|> 0$ and $|\widetilde{w}_0|> 0$ in $(0,1)$.
In fact, $J_{\nu_\theta}(s)>0$ for $0<s<j_{\nu_\theta,1}$. Moreover, from $|w_0|> 0$ in $(a_1,1)$ and 
$ |\widetilde{w}_0|> 0$ in $(a_2,1)$, we get 
$\begin{pmatrix}
    U_1^0(a_1)\\
    V_1^0(a_1)
\end{pmatrix}\neq \begin{pmatrix}
    0\\
    0
\end{pmatrix}$
and
$\begin{pmatrix}
    \widetilde{U}_1^0(a_2)\\
    \widetilde{V}_1^0(a_2)
\end{pmatrix}\neq \begin{pmatrix}
    0\\
    0
\end{pmatrix}$ in $(a_1,1)$ and $(a_2,1)$, respectively.
Using the hypothesis also on $(0,a_1)$ and $(0,a_2)$, we can conclude the same in $(0,a_1)$ and $(0,a_2)$.
 \end{Remark}

 \begin{Theorem}\label{theorem-strongly-right}
     Let $\theta \in [1,2)$, $0<a_1,a_2 <1$ and $0<t_1<t_2$. 
Let $w^{a_1}$ and $w^{a_2}$ be two solutions to \eqref{eq.deg-right-uniq}, corresponding to the initial values $w_0=u_0+iv_0$ and $\widetilde{w}_0=\widetilde{u}_0+i\widetilde{v}_0$, respectively. Assume that 
\begin{equation*}
\begin{pmatrix}
    U_1^0(a_1) \\
    V_1^0(a_1)
\end{pmatrix}\neq
  \begin{pmatrix}
    0 \\
    0
\end{pmatrix}
 \qquad \text{and} \qquad 
 \begin{pmatrix}
    \widetilde{U}_1^0(a_2) \\
    \widetilde{V}_1^0(a_2)
\end{pmatrix}\neq
  \begin{pmatrix}
    0 \\
    0
\end{pmatrix}, 
\end{equation*}
where $U^0_1(a_1)$ and $\widetilde{U}^0_1(a_2)$, for $u_0$ and $\widetilde{u}_0$ respectively, are given by~\eqref{eq.U0n} (with $n=1$) and $V^0_1(a_1)$ and $\widetilde{V}^0_1(a_2)$, for $v_0$ and $\widetilde{v}_0$ respectively, are given by~\eqref{eq.V0n} (with $n=1$).
 
Then $\partial_x w^{a_1}(1,t) = \partial_x w^{a_2}(1,t)$ for $t_1<t<t_2$ implies that $a_1=a_2$ and $u_0 = \widetilde{u}_0$, $v_0 = \widetilde{v}_0$ in $(a,1)$. 
 \end{Theorem}

 \textbf{Proof of Theorem~\ref{theorem-strongly-right}:} The proof is included in the proof of Theorem \ref{theorem2}.
 \hfill$\blacksquare$
 
\medskip 

We now state the uniqueness result for the third inverse problem, where the degeneracy power and the initial data are unknown.

\begin{Theorem}\label{theorem-strongly-right-theta}
     Let $\theta_1,\theta_2 \in [1,2)$ and $0<t_1<t_2$. 
Let $w^{\theta_1}$ and $w^{\theta_2}$ be two solutions to \eqref{eq.deg-right-uniq}, corresponding to the initial values $w_0=u_0+iv_0$ and $\widetilde{w}_0=\widetilde{u}_0+i\widetilde{v}_0$, respectively. Assume that 
\begin{equation}\label{eq-uniq-theta}
\begin{pmatrix}
    U_1^0(a;\theta_1) \\
    V_1^0(a;\theta_1)
\end{pmatrix}\neq
  \begin{pmatrix}
    0 \\
    0
\end{pmatrix}
 \qquad \text{and} \qquad 
 \begin{pmatrix}
    \widetilde{U}_1^0(a;\theta_2) \\
    \widetilde{V}_1^0(a;\theta_2)
\end{pmatrix}\neq
  \begin{pmatrix}
    0 \\
    0
\end{pmatrix}, 
\end{equation}
where $U^0_1(a;\theta_1)$ and $\widetilde{U}^0_1(a;\theta_2)$, for $u_0$ and $\widetilde{u}_0$ respectively, are given by~\eqref{eq.U0n} (with $n=1$) and $V^0_1(a;\theta_1)$ and $\widetilde{V}^0_1(a;\theta_2)$, for $v_0$ and $\widetilde{v}_0$ respectively, are given by~\eqref{eq.V0n} (with $n=1$).
 
Then $\partial_x w^{\theta_1}(1,t) = \partial_x w^{\theta_2}(1,t)$ for $t_1<t<t_2$ implies that $\theta_1=\theta_2$ and $u_0 = \widetilde{u}_0$, $v_0 = \widetilde{v}_0$ in $(a,1)$. 
 \end{Theorem}

Before we prove this uniqueness result, we need to prove the monotonicity of the first eigenvalue with respect to the degeneracy power $\theta$.

\subsection{Monotonicity of the first eigenvalue with respect to $\theta$}

In this subsection, we will consider the right eigenvalue problem \eqref{spectral-pb-right}
and prove the monotonicity of the first eigenvalue $\lambda_1$ w.r.t. $\theta$.
Since $\theta<2$, the injection of $H^1_\theta(a,1)$ into $L^2(a,1)$ is compact (see \cite{Al-Bou-etal} and \cite{PMC-Frr-Mrt-19} for the degeneracy coefficient $x^\theta$ in $(0,1)$), therefore, we can apply the spectral theorem in the right interval $(a,1)$ and give the variational formulation of the first eigenvalue:
\begin{equation*}
    \lambda_1(\theta):=\min_{\phi_1\in H^1_\theta\setminus \{0\}}\dfrac{\int_a^1 \left( x-a \right) ^\theta |\phi'_1|^2\, dx}{\int_a^1 |\phi_1|^2\, dx},
\end{equation*}
where $\phi_1$ is an associated eigenfunction.

First of all, we want to prove that if $\theta_1\leq\theta_2$, then $\lambda_1(\theta_1)\geq \lambda_1(\theta_2)$.
Observe that, for $x\in(a,1)$, $\theta_1\leq\theta_2$ implies $(x-a)^{\theta_2}\leq (x-a)^{\theta_1}$.
Hence, taking $\phi_1\in H^1_{\theta_1}(a,1)$, we get
\begin{equation*}
    \int_a^1 \left( x-a \right) ^{\theta_2} |\phi'_1|^2\, dx \leq \int_a^1 \left( x-a \right) ^{\theta_1} |\phi'_1|^2\, dx 
\end{equation*}
and
\begin{equation*}
    \|\phi_1\|^2_{H^1_{\theta_2}}=\|\phi_1\|^2_{L^2}+\int_a^1 \left( x-a \right) ^{\theta_2} |\phi'_1|^2\, dx \leq \|\phi_1\|^2_{L^2}+\int_a^1 \left( x-a \right) ^{\theta_1} |\phi'_1|^2\, dx
    = \|\phi_1\|^2_{H^1_{\theta_1}} < +\infty,
\end{equation*}
obtaining the continuous injection $H^1_{\theta_1}\subseteq  H^1_{\theta_2}$.
Thus, taking the minimum over $\phi_1\in H^1_{\theta_1}\setminus\{0\}$, we have
\begin{equation*}
    \min_{\phi_1\in H^1_{\theta_1}\setminus \{0\}}\dfrac{\int_a^1 \left( x-a \right) ^{\theta_2} |\phi'_1|^2 dx}{\int_a^1 |\phi_1|^2 dx} \leq \min_{\phi_1\in H^1_{\theta_1}\setminus \{0\}}\dfrac{\int_a^1 \left( x-a \right) ^{\theta_1} |\phi'_1|^2 dx}{\int_a^1 |\phi_1|^2 dx}:=\lambda_1(\theta_1).
\end{equation*}
We now want to prove that
\begin{equation}\label{eq-min}
    \lambda_1(\theta_2):=\min_{\phi_1\in H^1_{\theta_2}\setminus \{0\}}\dfrac{\int_a^1 \left( x-a \right) ^{\theta_2} |\phi'_1|^2\, dx}{\int_a^1 |\phi_1|^2\, dx}=\min_{\phi_1\in H^1_{\theta_1}\setminus \{0\}}\dfrac{\int_a^1 \left( x-a \right) ^{\theta_2} |\phi'_1|^2\, dx}{\int_a^1 |\phi_1|^2\, dx},
\end{equation}
observing that, due to the injection $H^1_{\theta_1}\subseteq  H^1_{\theta_2}$, we have
\begin{equation}\label{dis-lambda}
    \lambda_1(\theta_2)\leq \min_{\phi_1\in H^1_{\theta_1}\setminus \{0\}}\dfrac{\int_a^1 \left( x-a \right) ^{\theta_2} |\phi'_1|^2\, dx}{\int_a^1 |\phi_1|^2\, dx}.
\end{equation}
To obtain the equality, we exploit the density of $H^1_{\theta_1}$ in $H^1_{\theta_2}$.
In fact, the functions of class $C^\infty$ compactly supported in $(a,1)$ are dense in each $H^1_{\theta}$ (see \cite{Campiti-etal}).
Therefore, from $C^\infty_0(a,1)\subseteq H^1_{\theta_1}(a,1) \subseteq H^1_{\theta_2}(a,1)$, we deduce $\overline{H^1_{\theta_1}}^{\|\cdot\|_{H^1_{\theta_2}}}=H^1_{\theta_2}$.
Using this density result, 
we can also conclude that
\begin{equation*}
   \lambda_1(\theta_2)\geq \min_{\phi_1\in H^1_{\theta_1}\setminus\{0\}} \dfrac{\int_a^1 \left( x-a \right) ^{\theta_2} |\phi'_1|^2\, dx}{\int_a^1 |\phi_1|^2\, dx} 
\end{equation*}
and, using also \eqref{dis-lambda}, we deduce the validity of equality \eqref{eq-min}, from which we get $\lambda_1(\theta_2)\leq\lambda_1(\theta_1)$.

Finally, we want to prove the strict monotonicity of $\lambda_1$ w.r.t. $\theta$.
Taking into account a minimum point $\tilde{\phi}_1\in H^1_{\theta_1}\setminus\{0\}$ and $\theta_1<\theta_2$, we obtain $(x-a)^{\theta_2}< (x-a)^{\theta_1}$ for $x\in(a,1)$.
Hence, $\int_a^1 \left(\left( x-a \right) ^{\theta_1}-\left( x-a \right) ^{\theta_2}\right) |\tilde{\phi}^\prime_1|^2 dx=0$ would imply $|\tilde{\phi}^\prime_1|^2=0$ a.e. on $(a,1)$, which is a contradiction because it would also require $\tilde{\phi}_1=0$.
Thus, we obtain $\lambda_1(\theta_2)<\lambda_1(\theta_1)$.
\medskip 

 \textbf{Proof of Theorem~\ref{theorem-strongly-right-theta}:}
Let us consider the right sub-problem \eqref{eq.deg-right-uniq} and set 
\begin{equation*}\label{lambdamu-theta}
\lambda_n:= k_{\theta_1} ^2 \frac{j_{\nu _{\theta_1} ,n}^2}{(1-a)^{2k_{\theta_1}}}, \quad \mu_n:= k_{\theta_2} ^2 \frac{j_{\nu _{\theta_2} ,n}^2}{(1-a)^{2k_{\theta_2}}}, \quad n\in \mathbb{N},
\end{equation*}
for which we have $\lambda_1<\lambda_2<\ldots$ and $\mu_1<\mu_2<\ldots$. 
Assume that $\theta_1\neq \theta_2$. Without loss of generality, we can assume that $\theta_1<\theta_2$, then $\lambda_1>\mu_1$, due to the above monotonicity result.
 Using, for example, the condition $\partial_x w^{\theta_1}(1,t) = \partial_x w^{\theta_2}(1,t)$ and the analyticity w.r.t. time, we now prove that $\theta_1$=$\theta_2$.
 Therefore, setting
\begin{equation*}
     W^0_n(a;\theta_1):=\begin{pmatrix} 
    U_n^0(a;\theta_1)\\
    V_n^0(a;\theta_1)
\end{pmatrix}
\quad \text{and} \quad \widetilde{W}^0_n(a;\theta_2):=\begin{pmatrix}
    \widetilde{U}^0_n(a;\theta_2)\\
    \widetilde{V}^0_n(a;\theta_2)
\end{pmatrix},
\end{equation*}
we obtain 
\begin{equation*}\label{eq2.uniq2t2-theta}
\displaystyle\sum_{n = 1}^\infty \dfrac{ e^{-\lambda_n t}}{d_{\nu_{\theta_1},n}} R(\beta t)
W^0_n(a;\theta_1)
=\displaystyle\sum_{n = 1}^\infty \dfrac{ e^{-\mu_n t}}{d_{\nu_{\theta_2},n}} R(\tilde{\beta} t)
\widetilde{W}^0_n(a;\theta_2)
\end{equation*}
for $t>t_1$.
Therefore, taking the absolute value, we have 
\begin{equation*}
   \left| \dfrac{e^{-(\lambda_1-\mu_1) t}}{d_{\nu_{\theta_1},1}}
W^0_1(a;\theta_1)+\displaystyle\sum_{n = 2}^\infty \dfrac{e^{-(\lambda_n-\mu_1) t}}{d_{\nu_{\theta_1},n}}
W^0_n(a;\theta_1)\right|
=\left|\dfrac{1}{d_{\nu_{\theta_2},1}} 
\widetilde{W}^0_1(a;\theta_2)+\displaystyle\sum_{n = 2}^\infty \dfrac{e^{-(\mu_n-\mu_1) t}}{d_{\nu_{\theta_2},n}}
\widetilde{W}^0_n(a;\theta_2)\right|.
\end{equation*}
Since $\lambda_1>\mu_1$, we let $t\to\infty$ in the previous equality to get
$
\left|
    \widetilde{W}^0_1(a;\theta_2)\right|=
    0 $,
in contrast with \eqref{eq-uniq-theta}. Thus, $\theta_1=\theta_2$ follows, and this implies $\lambda_n=\mu_n$ for all $n\in \mathbb{N}$.
The uniqueness of the initial data in $(a,1)$ can be achieved following the same reasoning as in the Proof of Theorem \ref{theorem2}.

\begin{Remark}
    We observe that the above method can be adapted to obtain a uniqueness result for a more general class of problems, even in abstract form. For this purpose one could consider a family of differential operators such that the corresponding first eigenvalues satisfy a suitable monotonicity result, w.r.t. the unknown parameter that we are reconstructing, and a non degenerate observation of all first eigenfunctions exists. We will not further develop such an investigation in this paper, leaving it for a forthcoming study.
\end{Remark}

\section{Real systems of 1-D coupled degenerate parabolic equations}\label{sec:system}

The degenerate parabolic equation \eqref{eq.deg}, with complex solution $w(x,t)=u(x,t)+iv(x,t)$, can also be reformulated as a real system of 1-D coupled degenerate parabolic equations with the following structure: 
\begin{equation}\label{system}
\begin{cases}
\partial_t u- \partial_x(|x-a|^\theta \partial_x u) -\alpha u+\beta v=0, & (x,t)\in(0,1)\times(0,T),  \\[2mm]
\partial_t v- \partial_x(|x-a|^\theta \partial_x v) -\alpha v-\beta u=0, & (x,t)\in(0,1)\times(0,T),  \\[1mm]
\begin{pmatrix}
    u(0,t)\\
    v(0,t)
\end{pmatrix} =\begin{pmatrix}
    0\\
    0
\end{pmatrix} , \quad \begin{pmatrix}
    u(1,t)\\
    v(1,t)
\end{pmatrix} =\begin{pmatrix}
    0\\
    0
\end{pmatrix}, &  t\in (0,T),  \\[5mm]
\begin{pmatrix}
    u(x,0)\\
    v(x,0)
\end{pmatrix} =\begin{pmatrix}
    u_0(x)\\
    v_0(x)
\end{pmatrix},  & x\in(0,1).
\end{cases}
\end{equation}
Focusing now on the right interval $(a,1)$, we can analyze the sub-system
\begin{equation}\label{sub-system}
\begin{cases}
\partial_t u- \partial_x((x-a)^\theta \partial_x u) -\alpha u+\beta v=0, & (x,t)\in(a,1)\times(0,T),  \\[2mm]
\partial_t v- \partial_x((x-a)^\theta \partial_x v) -\alpha v-\beta u=0, & (x,t)\in(a,1)\times(0,T),  \\[1mm]
\begin{pmatrix}
  (x-a)^\theta \partial_x u|_{x=a}  \\
   (x-a)^\theta \partial_x v|_{x=a} 
\end{pmatrix}=\begin{pmatrix}
    0\\0
\end{pmatrix}, \quad\begin{pmatrix}
    u(1,t)\\
    v(1,t)
\end{pmatrix} =\begin{pmatrix}
    0\\
    0
\end{pmatrix}, &  t\in (0,T),  \\[4mm]
\begin{pmatrix}
    u(x,0)\\
    v(x,0)
\end{pmatrix} =\begin{pmatrix}
    u_0(x)\\
    v_0(x)
\end{pmatrix},  & x\in(a,1).
\end{cases}
\end{equation}
The solution of the system \eqref{sub-system} reads
    $\begin{pmatrix}
        u(x,t)\\
        v(x,t)
    \end{pmatrix}=e^{\alpha t}R(\beta t)
    \begin{pmatrix}
        e^{tA}u_0 \\
        e^{tA}v_0
    \end{pmatrix}$,
where $R(\beta t)$ is the rotation matrix defined in \eqref{rotation} and $e^{tA}u_0$, $e^{tA}v_0$ are defined by \eqref{eq.sol.u} and \eqref{eq.sol.v}.
The vector of normal derivatives at the boundary $\begin{pmatrix}
        \partial_x u(1,t)\\
        \partial_x v(1,t)
    \end{pmatrix}$
is determined exploiting the two components \eqref{eq.normderuxa}, \eqref{eq.normdervxa}.

The Theorem that allows us to achieve the Lipschitz stability estimate can be stated as follows:
\begin{Theorem}\label{th.stability-system}
Let $\theta \in [1,2)$ and assume that $u_0,v_0\in Lip([0,1])$. Let $\begin{pmatrix}
    u^{a_1}\\
    v^{a_1}
\end{pmatrix}$ and $\begin{pmatrix}
    u^{a_2}\\
    v^{a_2}
\end{pmatrix}$ be the solutions to \eqref{sub-system} corresponding to the degeneracy points $a_1$ and $a_2$, respectively. Assume that there exist $\delta>0$ and $[\tau,\gamma]\subset (0,1)$ such that 
\begin{equation*}\label{th2.assump-system}
 \left|\begin{pmatrix}
    U_1^0(a)\\
    V_1^0(a)
\end{pmatrix}\right|\ge \delta, \qquad \forall\, a\in [\tau,\gamma],
\end{equation*}
with $U_1^0(a)$ and $V_1^0(a)$ given by \eqref{eq.U0n} and \eqref{eq.V0n} for $\theta \in [1,2)$ and $n=1$. Then, there exist $t_0(u_0,v_0,\delta,L,\theta)>0$ and a constant $C>0$ such that the stability estimate 
  \begin{equation*}\label{eq.th-system} 
 |a_2-a_1| \le   C \left| \begin{pmatrix}
   \partial_x u^{a_2}(1,t)\\
   \partial_x v^{a_2}(1,t)
\end{pmatrix} -  \begin{pmatrix}
   \partial_x u^{a_1}(1,t)\\
   \partial_x v^{a_1}(1,t)
\end{pmatrix}\right|  
 \end{equation*}
 holds 
 \begin{itemize}
     \item for all $a_1,a_2\in [\tau,\gamma]$ and for all $t\in  [t_0, t_1]$ (with $t_1>t_0$), if $\lambda_1(\gamma)>\alpha$;
     \item for all $a_1,a_2\in [\tau,\gamma]$ and for all $t\geq t_0$, if $\lambda_1(\gamma)\leq\alpha$,
 \end{itemize}
 where $\lambda_1(\gamma)=
     k_\theta^2\dfrac{j_{\nu_\theta,1}^2}{(1-\gamma)^{2k_\theta}}$.
\end{Theorem}
In order to also state the uniqueness result of the solution, we take into account also the left sub-problem
\begin{equation}\label{sub-system-left}
\begin{cases}
\partial_t u- \partial_x((a-x)^\theta \partial_x u) -\alpha u+\beta v=0, & (x,t)\in(0,a)\times(0,T),  \\[2mm]
\partial_t v- \partial_x((a-x)^\theta \partial_x v) -\alpha v-\beta u=0, & (x,t)\in(0,a)\times(0,T),  \\[1mm]
\begin{pmatrix}
    u(0,t)\\
    v(0,t)
\end{pmatrix} =\begin{pmatrix}
    0\\
    0
\end{pmatrix}, \quad \begin{pmatrix}
  (a-x)^\theta \partial_x u|_{x=a}  \\
   (a-x)^\theta \partial_x v|_{x=a} 
\end{pmatrix}=\begin{pmatrix}
    0\\0
\end{pmatrix}, &  t\in (0,T),  \\[4mm]
\begin{pmatrix}
    u(x,0)\\
    v(x,0)
\end{pmatrix} =\begin{pmatrix}
    u_0(x)\\
    v_0(x)
\end{pmatrix},  & x\in(0,a).
\end{cases}
\end{equation}
The normal derivatives at the boundary 
    $\begin{pmatrix}
        \partial_x u(0,t)\\
        \partial_x v(0,t)
    \end{pmatrix}$
are determined using \eqref{normalder-left}. 
Once again, we can distinguish between three uniqueness Theorems.
\begin{Theorem}\label{theorem2-uniq-system}
Let $\theta \in [1,2)$, $0<a_1,a_2 <1$ and $0<t_1<t_2$. 
Let $\begin{pmatrix}
    u^{a_1}\\
    v^{a_1}
\end{pmatrix}$ and $\begin{pmatrix}
    u^{a_2}\\
    v^{a_2}
\end{pmatrix}$ be two solutions to \eqref{sub-system} and \eqref{sub-system-left}, corresponding to the initial values $\begin{pmatrix}
    u_0\\
    v_0
\end{pmatrix}$ and $\begin{pmatrix}
    \widetilde{u}_0\\
    \widetilde{v}_0
\end{pmatrix}$, respectively, and the coefficients $c=\alpha+i\beta$ and $c=\tilde{\alpha}+i\tilde{\beta}$, respectively. Assume that 
\begin{equation*}\label{th2.eq1-system}
\begin{pmatrix}
    U_1^0(a_1) \\
    V_1^0(a_1)
\end{pmatrix}\neq
  \begin{pmatrix}
    0 \\
    0
\end{pmatrix}
 \qquad \text{and} \qquad 
 \begin{pmatrix}
    \widetilde{U}_1^0(a_2) \\
    \widetilde{V}_1^0(a_2)
\end{pmatrix}\neq
  \begin{pmatrix}
    0 \\
    0
\end{pmatrix}, 
\end{equation*}
where $U^0_1(a_1)$ and $\widetilde{U}^0_1(a_2)$ for $u_0$ and $\widetilde{u}_0$, respectively, are given by
\begin{itemize}
    \item 
   \eqref{eq.U0n} or \eqref{eq.U0n-left} (with $n=1$) for the sub-problem \eqref{sub-system} or \eqref{sub-system-left}, respectively;
   \end{itemize}
and $V^0_1(a_1)$ and $\widetilde{V}^0_1(a_2)$ for $v_0$ and $\widetilde{v}_0$, respectively, are given by
\begin{itemize}
    \item \eqref{eq.V0n} or \eqref{eq.V0n-left} (with $n=1$) for the sub-problem \eqref{sub-system} or \eqref{sub-system-left}, respectively.
\end{itemize}
\noindent 
Then $\begin{pmatrix}
   \partial_x u^{a_1}(1,t)\\
   \partial_x v^{a_1}(1,t)
\end{pmatrix} = \begin{pmatrix}
   \partial_x u^{a_2}(1,t)\\
   \partial_x v^{a_2}(1,t)
\end{pmatrix}$ and $\begin{pmatrix}
   \partial_x u^{a_1}(0,t)\\
   \partial_x v^{a_1}(0,t)
\end{pmatrix} = \begin{pmatrix}
   \partial_x u^{a_2}(0,t)\\
   \partial_x v^{a_2}(0,t)
\end{pmatrix}$ for $t_1<t<t_2$ imply that $\alpha=\tilde{\alpha}$, $a_1=a_2$, $u_0 = \widetilde{u}_0$, $v_0 = \widetilde{v}_0$ and $\beta=\tilde{\beta}$. 
 \end{Theorem}

 \begin{Theorem}\label{theorem2-uniq-system-2}
Let $\theta \in [1,2)$, $0<a_1,a_2 <1$ and $0<t_1<t_2$. 
Let $\begin{pmatrix}
    u^{a_1}\\
    v^{a_1}
\end{pmatrix}$ and $\begin{pmatrix}
    u^{a_2}\\
    v^{a_2}
\end{pmatrix}$ be two solutions to \eqref{sub-system}, corresponding to the initial values $\begin{pmatrix}
    u_0\\
    v_0
\end{pmatrix}$ and $\begin{pmatrix}
    \widetilde{u}_0\\
    \widetilde{v}_0
\end{pmatrix}$, respectively. Assume that 
\begin{equation*}\label{th2.eq1-system}
\begin{pmatrix}
    U_1^0(a_1) \\
    V_1^0(a_1)
\end{pmatrix}\neq
  \begin{pmatrix}
    0 \\
    0
\end{pmatrix}
 \qquad \text{and} \qquad 
 \begin{pmatrix}
    \widetilde{U}_1^0(a_2) \\
    \widetilde{V}_1^0(a_2)
\end{pmatrix}\neq
  \begin{pmatrix}
    0 \\
    0
\end{pmatrix}, 
\end{equation*}
where $U^0_1(a_1)$ and $\widetilde{U}^0_1(a_2)$, for $u_0$ and $\widetilde{u}_0$, respectively, are given by \eqref{eq.U0n} (with $n=1$);
and $V^0_1(a_1)$ and $\widetilde{V}^0_1(a_2)$, for $v_0$ and $\widetilde{v}_0$, respectively, are given by \eqref{eq.V0n} (with $n=1$).
\noindent 
Then $\begin{pmatrix}
   \partial_x u^{a_1}(1,t)\\
   \partial_x v^{a_1}(1,t)
\end{pmatrix} = \begin{pmatrix}
   \partial_x u^{a_2}(1,t)\\
   \partial_x v^{a_2}(1,t)
\end{pmatrix}$ for $t_1<t<t_2$ implies that $a_1=a_2$ and $u_0 = \widetilde{u}_0$, $v_0 = \widetilde{v}_0$ in $(a,1)$. 
 \end{Theorem}

 \begin{Theorem}\label{theorem2-uniq-system-theta}
Let $\theta_1, \theta_2 \in [1,2)$ and $0<t_1<t_2$. 
Let $\begin{pmatrix}
    u^{\theta_1}\\
    v^{\theta_1}
\end{pmatrix}$ and $\begin{pmatrix}
    u^{\theta_2}\\
    v^{\theta_2}
\end{pmatrix}$ be two solutions to \eqref{sub-system}, corresponding to the initial values $\begin{pmatrix}
    u_0\\
    v_0
\end{pmatrix}$ and $\begin{pmatrix}
    \widetilde{u}_0\\
    \widetilde{v}_0
\end{pmatrix}$, respectively. Assume that 
\begin{equation*}
\begin{pmatrix}
    U_1^0(a; \theta_1) \\
    V_1^0(a; \theta_1)
\end{pmatrix}\neq
  \begin{pmatrix}
    0 \\
    0
\end{pmatrix}
 \qquad \text{and} \qquad 
 \begin{pmatrix}
    \widetilde{U}_1^0(a; \theta_2) \\
    \widetilde{V}_1^0(a; \theta_2)
\end{pmatrix}\neq
  \begin{pmatrix}
    0 \\
    0
\end{pmatrix}, 
\end{equation*}
where $U^0_1(a; \theta_1)$ and $\widetilde{U}^0_1(a; \theta_2)$, for $u_0$ and $\widetilde{u}_0$, respectively, are given by \eqref{eq.U0n} (with $n=1$);
and $V^0_1(a; \theta_1)$ and $\widetilde{V}^0_1(a;\theta_2)$, for $v_0$ and $\widetilde{v}_0$, respectively, are given by \eqref{eq.V0n} (with $n=1$).
\noindent 
Then $\begin{pmatrix}
   \partial_x u^{\theta_1}(1,t)\\
   \partial_x v^{\theta_1}(1,t)
\end{pmatrix} = \begin{pmatrix}
   \partial_x u^{\theta_2}(1,t)\\
   \partial_x v^{\theta_2}(1,t)
\end{pmatrix}$ for $t_1<t<t_2$ implies that $\theta_1=\theta_2$ and $u_0 = \widetilde{u}_0$, $v_0 = \widetilde{v}_0$ in $(a,1)$. 
 \end{Theorem}

\section{Numerical results}\label{sec.numerics}

In this section, we will show some numerical results related to the identification of the degeneracy point $a\in(0,1)$, the initial data $(u_0,v_0)$ and also the degeneracy power $\theta$ in 
\begin{equation}\label{eq.cp} 
\begin{cases}
 \partial_t u - \partial_x(|x-a|^\theta \partial_x u) - \alpha u + \beta  v  = 0,  & (x,t)\in  (0,1)\times(0,T),   \\
 \partial_t v - \partial_x(|x-a|^\theta \partial_x v) - \alpha v - \beta  u  = 0,   & (x,t)\in  (0,1)\times(0,T),   \\ 
 u(0,t) = 0,  \quad u(1,t) = 0,    & t\in  (0,T),  \\
 v(0,t) = 0, \quad v(1,t) = 0,      & t\in  (0,T), \\
 u(x,0) = u_0(x),                       & x\in   (0,1), \\
 v(x,0) = v_0(x),                        & x\in  (0,1). \\
\end{cases}
\end{equation}
We will perform some numerical tests for the strong degeneracy case with $\theta\in [1,2)$, which have not been considered in the previous article \cite{PMC-Dou-Yam-24}, where only the scalar case with $\theta=1$ has been treated.
To make the graphical representation clearer, we have performed simulations with the diffusion term multiplied by a constant factor less than one.
In particular, this Section will be devoted to the numerical reconstruction of the solution for several kinds of inverse problem.
More precisely, the following two tests will consider the inverse problem of recovering the degeneracy point, for both one-point and distributed measurements. 

\begin{description}
\item[Test~1:] Find $a$ from the punctual measurements $\eta (t^\ast) = \partial_x u(1,t^\ast)$ and $\zeta (t^\ast) = \partial_x v(1,t^\ast)$ for some $t^\ast\in (0,T)$. 

\item [Test~2:] Find $a$ from distributed measurements  $\eta (t) =\partial_x u(1,t)$  for $t\in (t_1,t_2)$. 
\end{description}

The following three tests will consider the more general inverse problem of degeneracy and initial data reconstruction, taking into account distributed measurements on one or two sides of the domain.
These tests, especially Test $3$ and Test $5$, also illustrate the first two uniqueness theoretical results of Section~\ref{sec:uniqueness} and Section~\ref{sec:system}. 

\begin{description}
\item [Test~3:] Find $a$ and piecewise-constant initial data $(u_0,v_0)$ in $(0,1)$ from the distributed measurements $\eta (t) =\partial_x u(1,t)$ and $\zeta (t) =\partial_x v(1,t)$, $\rho (t) =\partial_x u(0,t)$ and $\kappa (t) =\partial_x v(0,t)$ for $t\in (t_1,t_2)$.

\item [Test~4:] Find $a$ and piecewise-constant initial data $(u_0,v_0)$ in $(0,1)$ from distributed measurements $\eta (t) =\partial_x u(1,t)$ and $\zeta (t) =\partial_x v(1,t)$ for $t\in (t_1,t_2)$. 

\item [Test~5:] Find $a$ and initial data $(u_0,v_0)$ in $(a,1)$ from the distributed measurements $\eta (t) =\partial_x u(1,t)$ and $\zeta (t) =\partial_x v(1,t)$ for $t\in (t_1,t_2)$.
\end{description}

The last test concerns the inverse problem of reconstructing the degeneracy power $\theta$ and the initial data from distributed measurements and covers the last uniqueness result in Section~\ref{sec:uniqueness} and Section~\ref{sec:system}.

\begin{description}
\item [Test~6:] Find $\theta$ and initial data $(u_0,v_0)$ in $(a,1)$ from the distributed measurements $\eta (t) =\partial_x u(1,t)$ and $\zeta (t) =\partial_x v(1,t)$ for $t\in (t_1,t_2)$.

\end{description}

\subsection{Degeneracy reconstruction with one-point measurements }

Given $T>0$ and $\eta (t^*)$ and $\zeta(t^*)$, we will present some numerical tests for the given initial data $u_0, v_0$, so as we can find $a\in (0,1)$ such that the solution to \eqref{eq.cp} for some $t^\ast\in (0,T)$ satisfies  
\begin{equation*}
\partial_x u(1,t^\ast) = \eta(t^\ast), 	\quad \partial_x v(1,t^\ast) = \zeta(t^\ast).
\end{equation*}

In order to reconstruct $a$, we will reformulate the inverse problem as an optimization problem. With fixed small $\delta>0$, let us consider the admissible set 

\begin{equation}\label{uad}
\mathcal{U}^{a}_{ad} =\{ a: a\in (\delta, 1-\delta) \} 
\end{equation}
and a functional $J : a\in \mathcal{U}^{a}_{ad} \mapsto \mathds{R}$ given by
$J(a) = \dfrac{1}{2}\displaystyle |\eta(t^\ast) -\partial_x u^a(1,t^\ast)|^2
+ \dfrac{1}{2}\displaystyle |\zeta(t^\ast) - \partial_x v^a(1,t^\ast)|^2$
for some $t^\ast\in (0,T)$.
\noindent
The related optimization problem is the following:
\begin{equation}\label{optpb1a}
\left\{\begin{array}{l}
\text{Minimize  $J(a) $, }   \\[1mm]
\text{where $a\in \mathcal{U}^{a}_{ad} $ and $(u^a,v^a)$ satisfies \eqref{eq.cp}. }
\end{array}\right.
\end{equation}
The \texttt{fmincon} function from MATLAB Optimization ToolBox (with a gradient based algorithm) will be used to solve the constrained optimization problem~\eqref{optpb1a}. 
The gradient of the cost functional is computed via the adjoint method. Our numerical implementation involves solving a system of four coupled equations to determine both the state variables and the derivatives with respect to the degeneracy point. We employ the MATLAB \texttt{pdepe} solver to handle the underlying partial differential equations. 

\begin{figure}[h!]
\begin{minipage}[t]{0.49\linewidth}
\noindent
\includegraphics[width=0.98\linewidth]{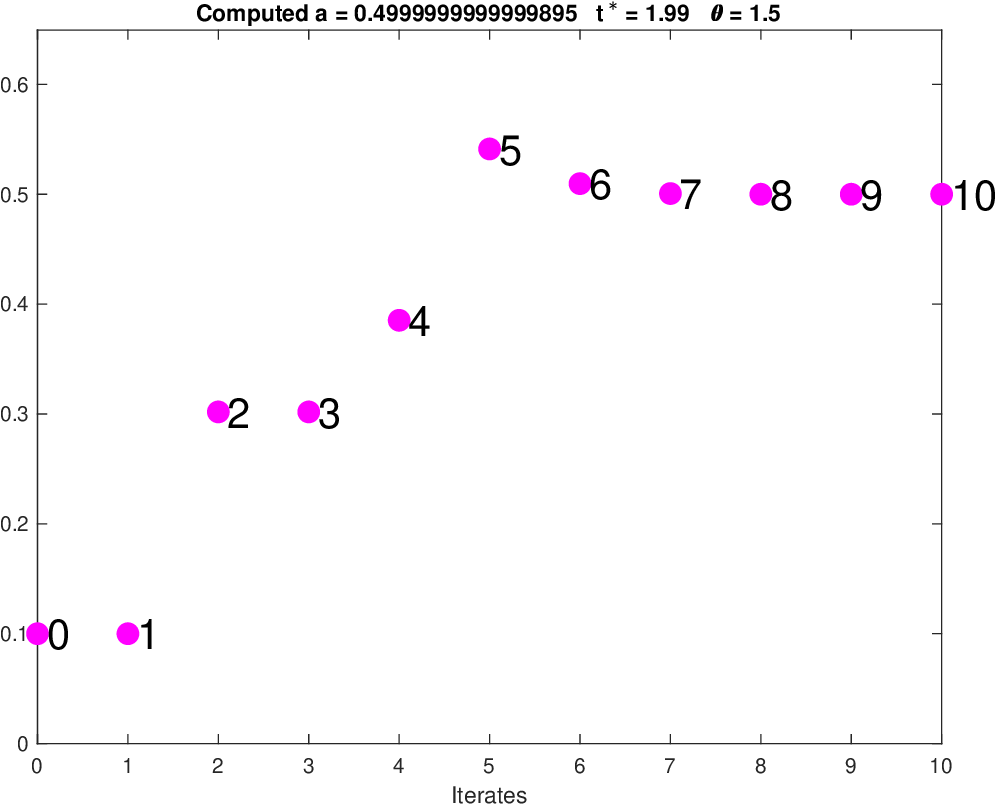}
\caption{Test~\ref{test1}, $\theta = 1.5$, $t^\ast = 1.99$, $u_0=1$, $v_0=1$. Iterations in the computation of $a$ by \texttt{trust-region-reflective} algorithm, $aini=0.1$. }
\label{fig.aTest1}
\end{minipage}
\hfill
\begin{minipage}[t]{0.49\linewidth}
\noindent
\includegraphics[width=0.98\linewidth]{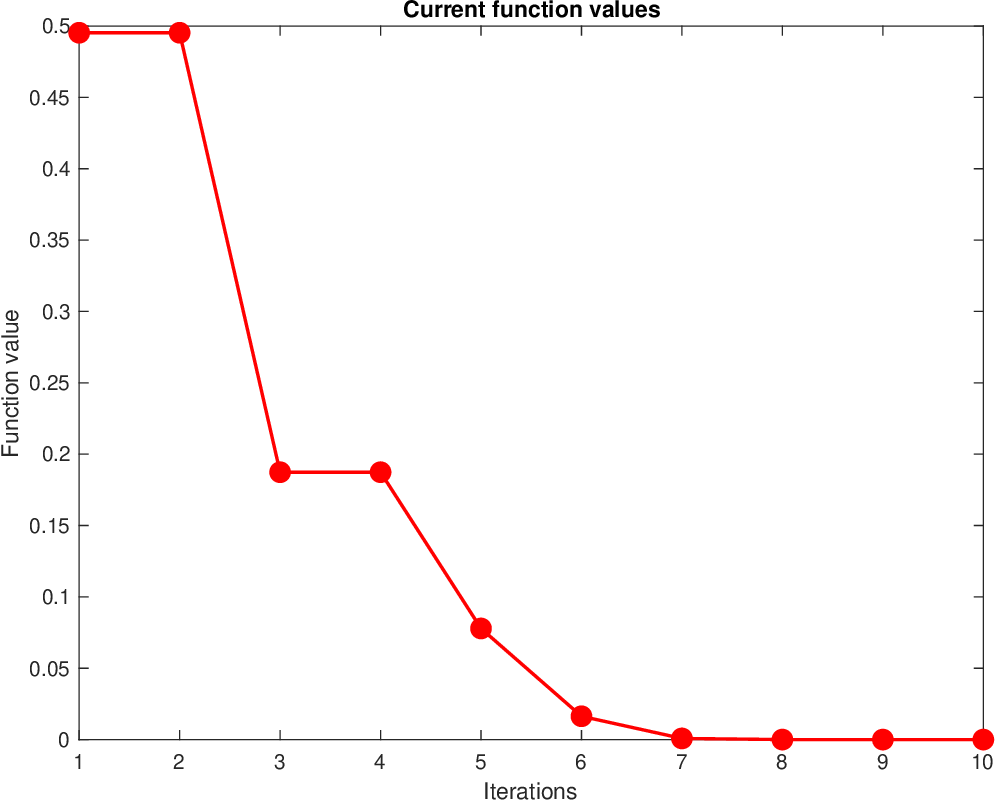}
\caption{Test~\ref{test1}, $\theta = 1.5$, $t^\ast = 1.99$, $u_0=1$, $v_0=1$. Evolution of the cost in \texttt{trust-region-reflective} algorithm, $aini=0.1$.}
\label{fig.fvalTest1}
\end{minipage}
\end{figure}

\begin{test} \label{test1}
The goal is to reconstruct the degeneracy point $a$ for a strong degenerate case and with the given initial data $u_0=1$, $v_0 =1$. 
We will take $\theta = 1.5$, $\alpha=1$, $\beta=1$, $T=4$, $t^\ast = 1.99$ and $aini=0.1$ as initial guess to recover the desired value of $a_d=0.5$ by the minimization algorithm.  
The numerical results can be seen in Figures~\ref{fig.aTest1} and~\ref{fig.fvalTest1}. The round points correspond to iterations during the optimization algorithm and the digits show the number of iterations performed. 
With the solid line, we have represented the evolution of the cost.  We obtain the computed value $a_c = 0.4999999999999895$ and the cost $J(a_c) \approx 1.e-27$. 

\end{test}
In Table~\ref{Table1} we can see the evolution of the cost when we introduce random uniform noises in the target. These results correspond to the \texttt{trust-region-reflective} algorithm. Notice that this method is highly optimized for problems where the gradient is provided via the adjoint method and it offers superior memory efficiency compared to alternative methods.

%

\begin{table}[h!]
\centering
\renewcommand{\arraystretch}{1.2}
{\begin{tabular}{cccc} \hline
Noise  & Cost & Iterations &  $a_c$ \\
\hline\hline
1\%         &  1.e-15          &    16      & 0.4967806209438190\\ \hline
0.1\%      &  1.e-16            &   10       &0.5010448098047946\\ \hline
0.01\%     &  1.e-19          &   10     & 0.5000110001604564  \\\hline
0\%          &   1.e-27         &    10      &  0.4999999999999895   \\  \hline
\end{tabular}}
\caption{Evolution of the cost with random noises in the target, Test~\ref{test1} with $\theta = 1.5$ and $aini=0.1$.}
\label{Table1}
\end{table}

For simplicity, we will not formally verify the numerical stability for the random noises for all other tests.
Whether theoretical stability can be demonstrated for all the above cases, potentially through alternative analytical techniques, has still to be proved.
\subsection{Degeneracy reconstruction with distributed measurements }

\begin{figure}[h!]
\begin{minipage}[t]{0.49\linewidth}
\noindent
\includegraphics[width=0.98\linewidth]{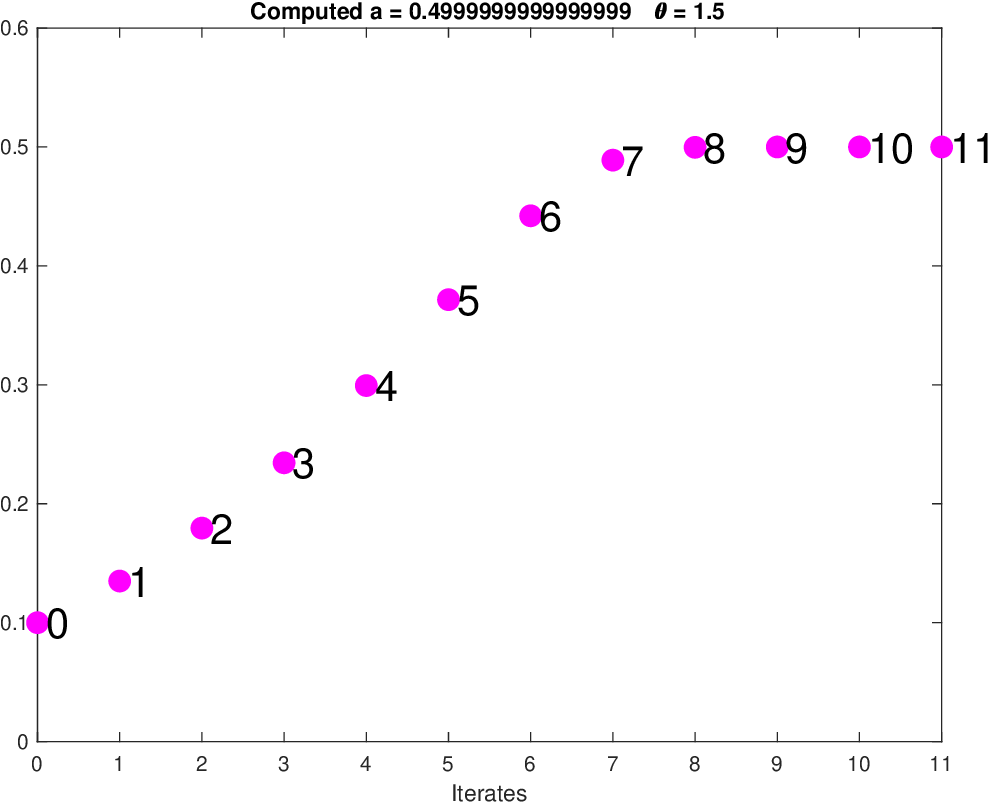}
\caption{Test~\ref{test2}, $\theta = 1.5$, $u_0=1$, $v_0=1$. Iterations in the computation of $a$ by \texttt{trust-region-reflective} algorithm, $aini=0.1$. }
\label{fig.aTest2}
\end{minipage}
\hfill
\begin{minipage}[t]{0.49\linewidth}
\noindent
\includegraphics[width=0.98\linewidth]{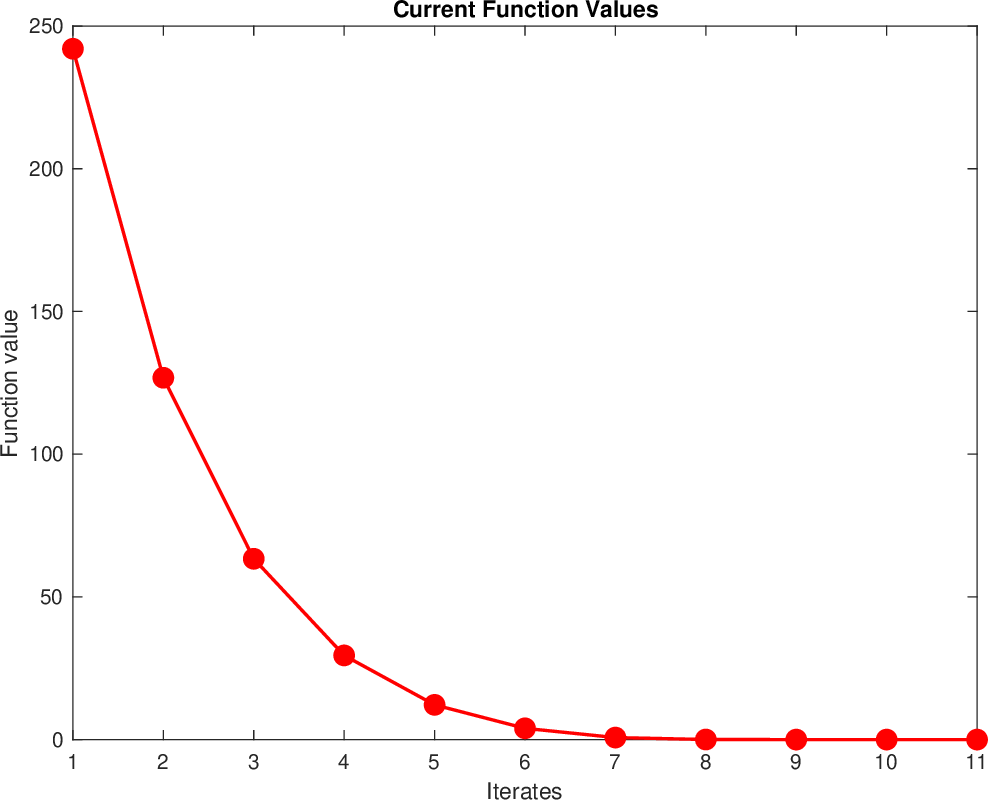}
\caption{Test~\ref{test2}, $\theta = 1.5$, $u_0=1$, $v_0=1$. Evolution of the cost in \texttt{trust-region-reflective} algorithm, $aini=0.1$.}
\label{fig.fvalTest2}
\end{minipage}
\end{figure}

 In this section, we will give some numerical simulations of reconstruction of the degeneracy point $a\in (0,1)$ only with one distributed measurement, so that, given $\eta(t)$, a solution to \eqref{eq.cp} satisfies
\begin{equation*}\label{obs1xo}
\partial_x u(1,t) = \eta(t), 	\quad \text{ for } t \in (t_1,t_2), \quad 0\le t_1 \le t_2 \le T.
\end{equation*}
As before, we reformulate the inverse problem as an optimization problem: 
\begin{equation*}
\left\{\begin{array}{l}
\text{Minimize  $\mathcal{I}(a) $, }   \\[1mm]
\text{where $a\in \mathcal{U}^{a}_{ad} $ and $(u^a,v^a)$ satisfies \eqref{eq.cp},}
\end{array}\right.
\end{equation*}
\noindent
where $\mathcal{U}^{a}_{ad}$ is given by \eqref{uad} and $\mathcal{I} : a\in \mathcal{U}^{a}_{ad} \mapsto \mathds{R}$ is defined as $\mathcal{I}(a) = \dfrac{1}{2}\displaystyle \int_{t_1}^{t_2} |\eta(t) -\partial_x u^a(1,t)|^2\, dt$.
\noindent
\begin{test} \label{test2}
We will take $\theta = 1.5$, $\alpha=1$, $\beta=1$, $T=4$, $t_1 = 0$, $t_2= T$, $u_0=1$, $v_0=1$ and $aini=0.1$ as initial guess to recover the desired value of $a_d=0.5$ using the minimization algorithm.  
The numerical results can be seen in Figures~\ref{fig.aTest2} and~\ref{fig.fvalTest2}. The round points correspond to iterations during the optimization algorithm.  With the solid line, we have represented the evolution of the cost.  We obtain the computed value $a_c = 0.4999999999999999$ and the cost $\mathcal{I}(a_c) \approx 1.e-26$.
\end{test}

\begin{Remark}
    We observe that in this simulation, the reconstruction of the degeneracy point is performed by measuring a single component of the normal derivative. This is always possible when $\beta$ is non-zero, using the same proof as the uniqueness Theorem \ref{theorem2}, but adapted for a single component. In the case where $\beta$ is zero, the two equations of the system are uncoupled, and reconstruction with a single component is possible only if at least one of the two vector components in the hypothesis \eqref{th2.eq1} is non-zero throughout the entire interval $[\tau,\gamma]$.
\end{Remark}
\subsection{Degeneracy and initial data reconstruction with distributed measurements }

 In this section, we will give some numerical simulations of reconstruction of the degeneracy point $a\in (0,1)$ and initial data $(u_0,v_0)$ by distributed measurements.
 The following analysis will allow us to better highlight the difference between the first two uniqueness Theorems in Section~\ref{sec:uniqueness} (or Section~\ref{sec:system} for the system version), in terms of the results for the initial data. 
 

 \subsubsection{Piecewise-constant initial data, two side measurements }

\begin{figure}[h!]
\begin{minipage}[t]{0.49\linewidth}
\noindent
\includegraphics[width=0.98\linewidth]{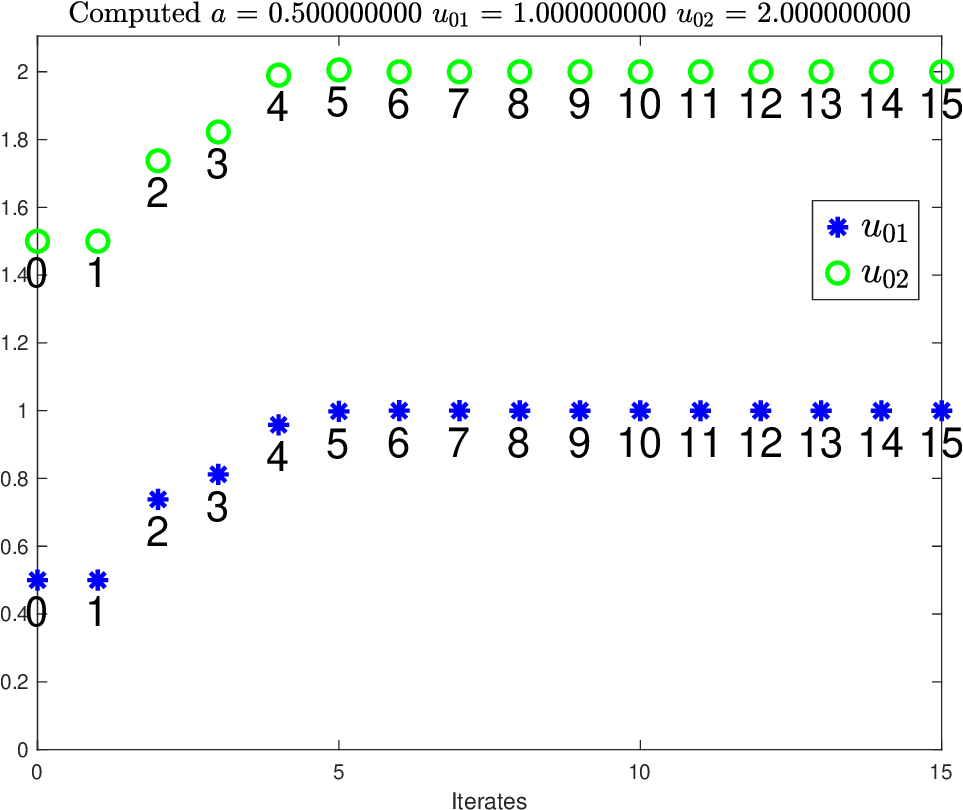}
\caption{Test~\ref{test4}, $\theta = 1.5$. Iterations in the computation of $u_{01}$ and $u_{02}$ by \texttt{trust-region-reflective} algorithm, $aini=0.1$, $u_{01ini}= 0.5$, $u_{02ini}= 1.5$. }
\label{fig.aTest4}
\end{minipage}
\hfill
\begin{minipage}[t]{0.49\linewidth}
\noindent
\includegraphics[width=0.98\linewidth]{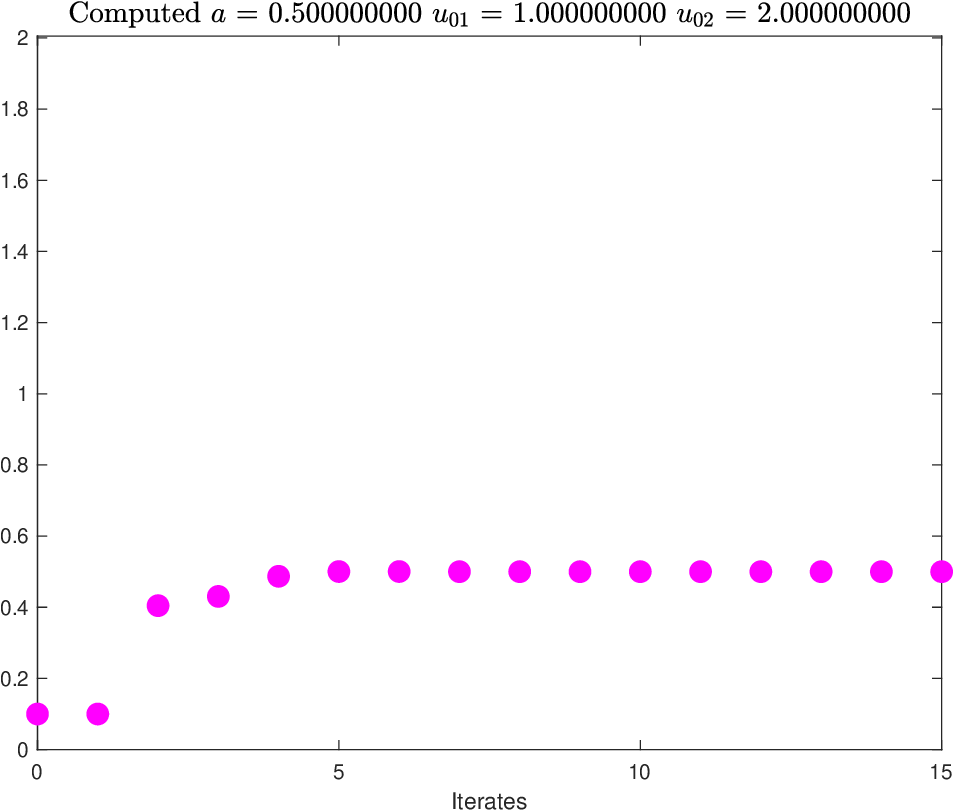}
\caption{Test~\ref{test4}, $\theta = 1.5$.  Iterations in the computation of $a$ by \texttt{trust-region-reflective} algorithm, $aini=0.1$, $u_{01ini}= 0.5$, $u_{02ini}= 1.5$.}
\label{fig.aTest4BB}
\end{minipage}
\end{figure}

\begin{figure}[h!]
\centering
\includegraphics[width=0.49\linewidth]{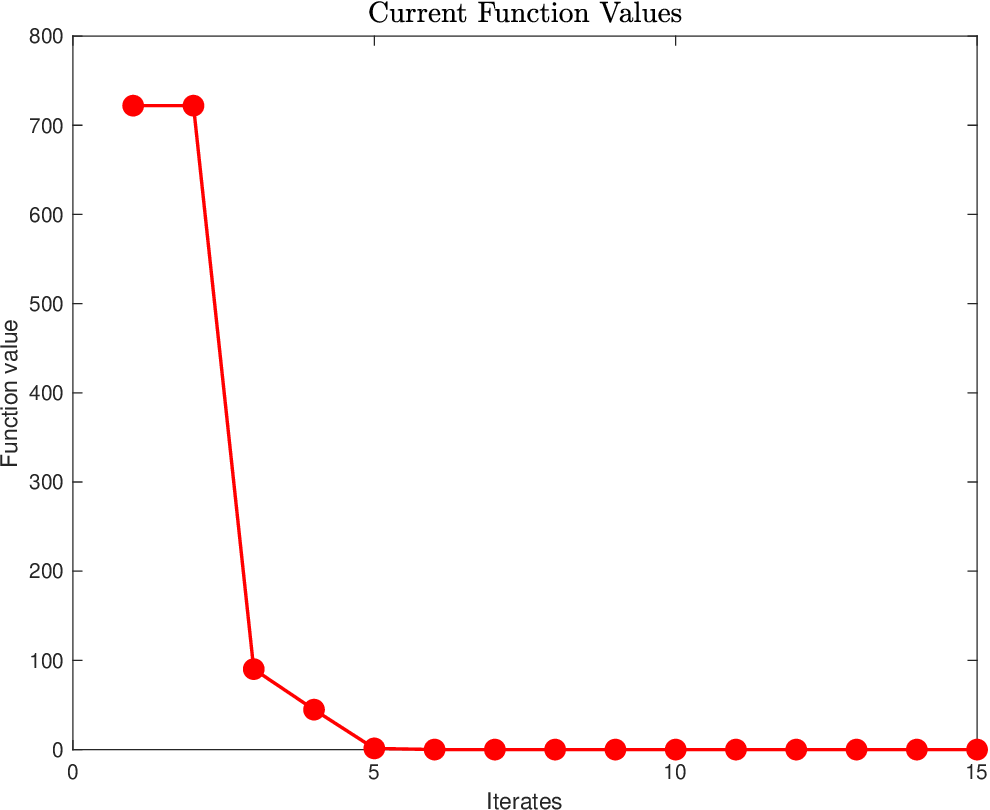}
\caption{Test~\ref{test4}, $\theta = 1.5$. Evolution of the cost in \texttt{trust-region-reflective} algorithm, $aini=0.1$, $u_{01ini}= 0.5$, $u_{02ini}= 1.5$.}
\label{fig.fvalTest4}%
\end{figure}

Let us present here the case where we can have some discontinuity in the initial data at the degeneracy point $a$. More precisely, we will assume that $u_0$ is of the form 
\begin{equation}
u_0 = \begin{cases} \label{inidata}     
u_{01}  & \text{if }   0< x<a, 
\\
u_{02}  & \text{if }  a<x<1,
\end{cases}
\end{equation}
with $u_{01}$ and $u_{02}$ constant
and, for simplicity, $v_0=0$. Therefore, our goal is to find $a\in (0,1)$ and initial data $u_{01}$ and $u_{02}$ such that a solution to \eqref{eq.cp} satisfies
\begin{equation*}
\begin{cases}
\partial_x u(0,t) = \rho(t), 	\quad  \quad \partial_x v(0,t) = \kappa(t),\\
\partial_x u(1,t) = \eta(t), 	\quad  \quad \partial_x v(1,t) = \zeta(t), 
\end{cases}
\quad \text{ for } t \in (t_1,t_2), \quad 0\le t_1 \le t_2 \le T. 
\end{equation*}
The reformulation of the inverse problem is as follows: 
\begin{equation*}
\left\{\begin{array}{l}
\text{Minimize  $\mathcal{H}(a,u_{01},u_{02})$, }  \\[1mm]
\text{where $a\in \mathcal{U}^{a}_{ad} $ and $(u^{a,u_{01},u_{02}},v^{a,u_{01},u_{02}})$ satisfies \eqref{eq.cp},}
\end{array}\right.
\end{equation*}
\noindent
where $\mathcal{U}^{a}_{ad}$ is given by \eqref{uad} and $\mathcal{H} :(a,u_{01},u_{02})\in \mathcal{U}^{a}_{ad}\times \mathds{R}\times\mathds{R} \mapsto \mathds{R}$ is defined as follows:
\begin{equation*}
\begin{array}{ll}
\mathcal{H}(a,u_{01},u_{02}) 
& =  \dfrac{1}{2}\displaystyle \int_{t_1}^{t_2} |\rho(t) -\partial_x u^{a,u_{01},u_{02}}(0,t)|^2\, dt +
 \dfrac{1}{2}\displaystyle \int_{t_1}^{t_2} |\kappa(t) -\partial_x v^{a,u_{01},u_{02}}(0,t)|^2\, dt 
 \\[5mm]
& + \dfrac{1}{2}\displaystyle \int_{t_1}^{t_2} |\eta(t) -\partial_x u^{a,u_{01},u_{02}}(1,t)|^2\, dt
   +
 \dfrac{1}{2}\displaystyle \int_{t_1}^{t_2} |\zeta(t) -\partial_x v^{a,u_{01},u_{02}}(1,t)|^2\, dt.
 \end{array}
\end{equation*}
\noindent
\begin{test} \label{test4}
We will take $\theta = 1.5$, $\alpha=1$, $\beta=1$, $T=4$, $t_1 = 0$, $t_2= T$, $u_{01ini} = 0.5$, $u_{02ini}= 1.5$, and $aini=0.1$ as initial guesses to recover the desired value of $a_d=0.5$, $u_{01d}= 1$, $u_{02d}= 2$ using the minimization algorithm.  
The numerical results can be seen in Figures~\ref{fig.aTest4}, ~\ref{fig.aTest4BB} and~\ref{fig.fvalTest4}. In Figure~\ref{fig.aTest4}, the stars and the round points correspond to the iterations during the optimization algorithm in the computation of $u_{01}$ and $u_{02}$, respectively. 
In Figure~\ref{fig.aTest4BB}, the round points represent the iterations during the optimization algorithm in the computation of $a$.
With the solid line in Figure~\ref{fig.fvalTest4}, we have represented the evolution of the cost during iterations.  We obtain the computed value $a_c = 0.5000000000049045$, $u_{01c} = 1.0000000000141918$, $u_{02c} = 1.9999999999341966$  and the cost $\mathcal{H}(a_c, u_{01c},u_{02c}) \approx 1.e-18$. This test allows us to numerically get the uniqueness result of Theorem \ref{theorem2}.
\end{test}

 \subsubsection{Piecewise-constant initial data, one side measurements }

\begin{figure}[h!]
\begin{minipage}[t]{0.49\linewidth}
\noindent
\includegraphics[width=0.98\linewidth]{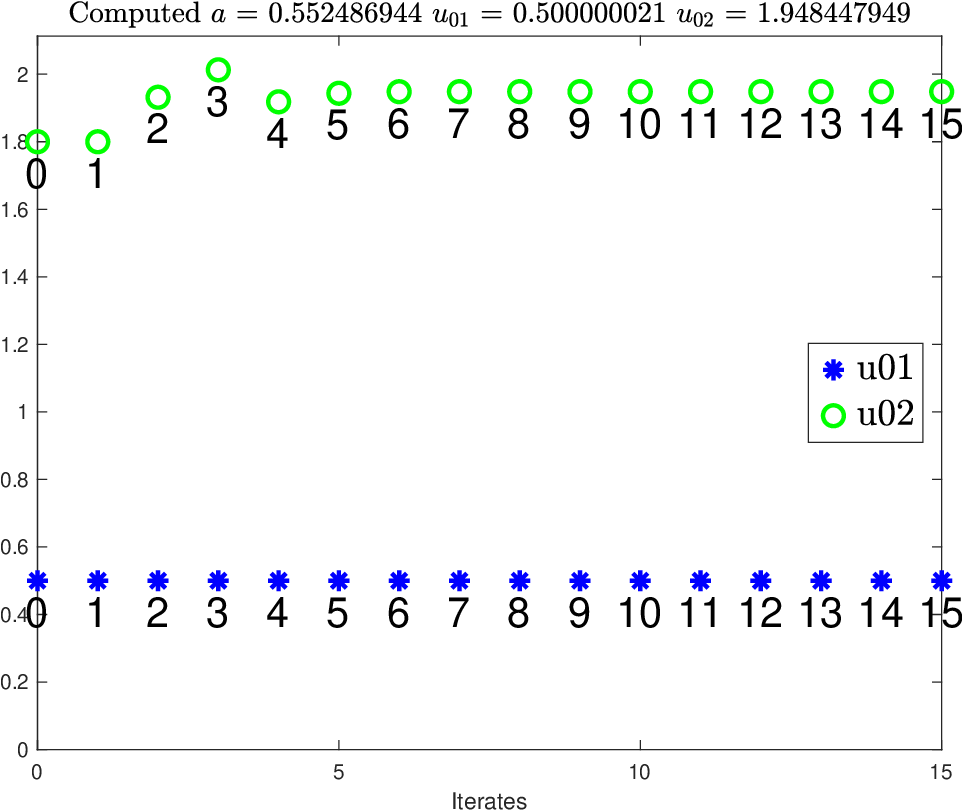}
\caption{Test~\ref{test5}, $\theta = 1.5$. Iterations in the computation of $u_{01}$ and $u_{02}$ by \texttt{trust-region-reflective} algorithm, $aini=0.1$, $u_{01ini}= 0.5$, $u_{02ini}= 1.8$. }
\label{fig.aTest5}
\end{minipage}
\hfill
\begin{minipage}[t]{0.49\linewidth}
\noindent
\includegraphics[width=0.98\linewidth]{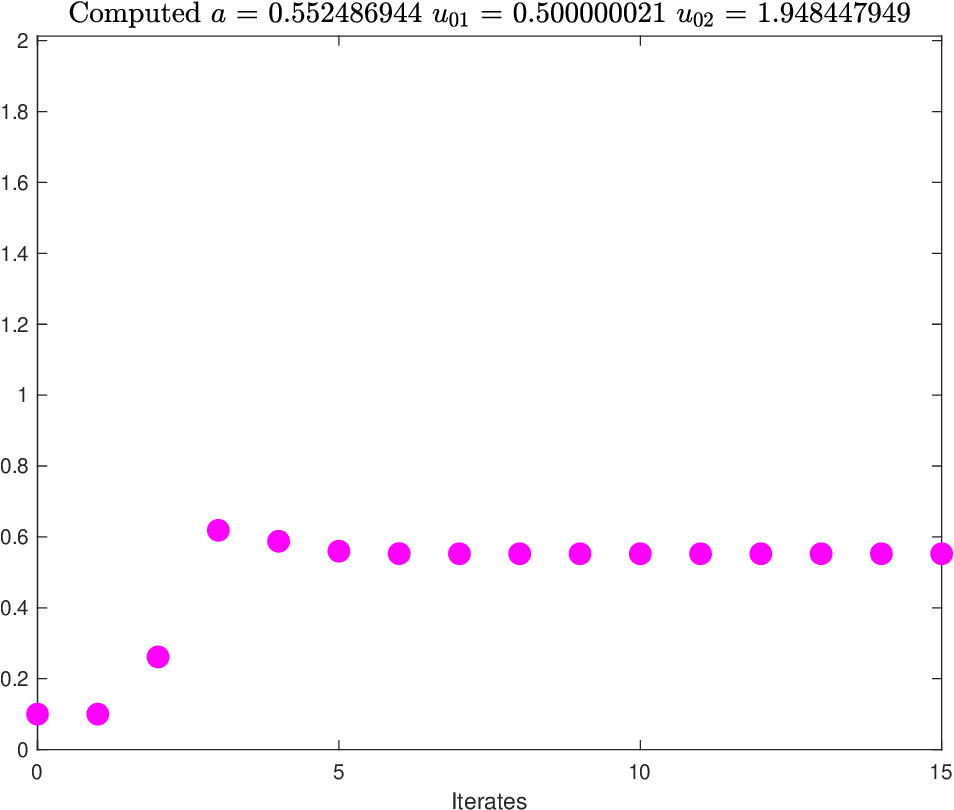}
\caption{Test~\ref{test5}, $\theta = 1.5$. Iterations in the computation of $a$ in \texttt{trust-region-reflective} algorithm, $aini=0.1$, $u_{01ini}= 0.5$, $u_{02ini}= 1.8$.}
\label{fig.aTest5BB}
\end{minipage}
\end{figure}

\begin{figure}[h!]
\centering
\includegraphics[width=0.49\linewidth]{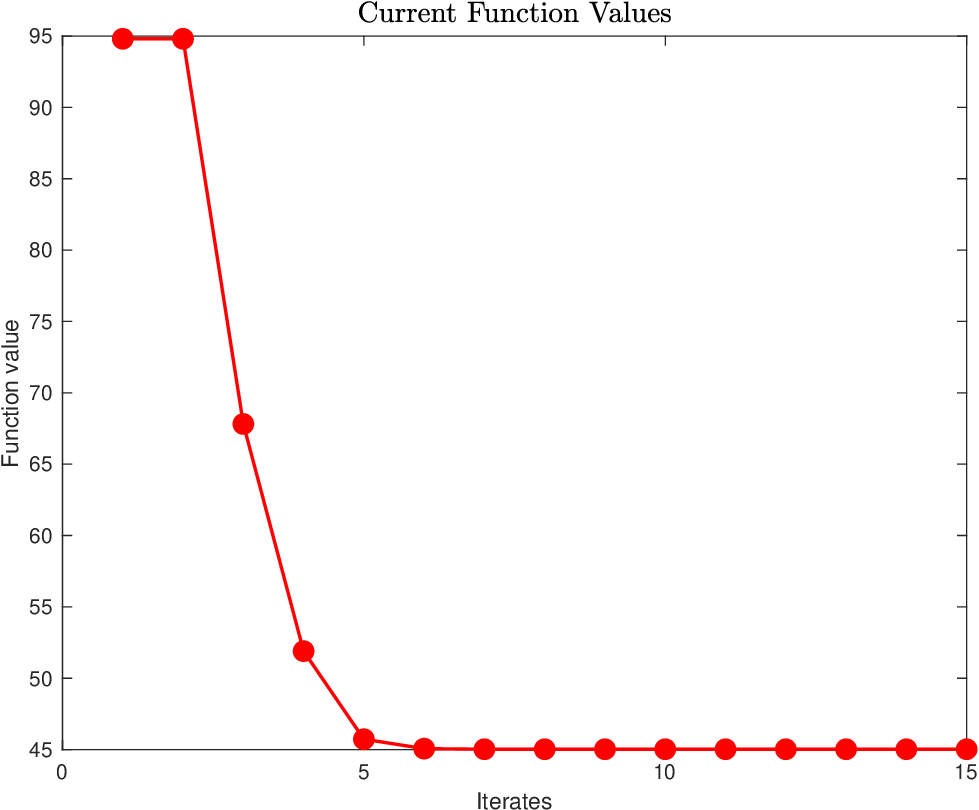}
\caption{Test~\ref{test5}, $\theta = 1.5$. Evolution of the cost in \texttt{trust-region-reflective} algorithm, $aini=0.1$, $u_{01ini}= 0.5$, $u_{02ini}= 1.8$.}
\label{fig.fvalTest5}
\end{figure}

Continuing in the case of discontinuous initial data with different constant on the two sides of the degeneracy point $a$, let us justify here that, for the reconstruction of the initial data on the whole interval $(0,1)$, one side measurement is not enough. 
More precisely, we will assume that $u_0$ is of the form \eqref{inidata} and, again $v_0=0$. Therefore, our goal is to find $a\in (0,1)$ and initial data $u_{01}$ and $u_{02}$ such that a solution to \eqref{eq.cp} satisfies
\begin{equation*}
\partial_x u(1,t) = \eta(t), 	\quad  \quad \partial_x v(1,t) = \zeta(t), 
\quad \text{ for } t \in (t_1,t_2), \quad 0\le t_1 \le t_2 \le T. 
\end{equation*}
The reformulation of the inverse problem is as follows: 
\begin{equation*}
\left\{\begin{array}{l}
\text{Minimize  $\mathcal{M}(a,u_{01},u_{02})$, }  \\[1mm]
\text{where $a\in \mathcal{U}^{a}_{ad} $ and $(u^{a,u_{01},u_{02}},v^{a,u_{01},u_{02}})$ satisfies \eqref{eq.cp},}
\end{array}\right.
\end{equation*}
\noindent
where $\mathcal{U}^{a}_{ad}$ is given by \eqref{uad} and $\mathcal{M} :(a,u_{01},u_{02})\in \mathcal{U}^{a}_{ad}\times \mathds{R}\times\mathds{R} \mapsto \mathds{R}$ is defined as follows:
\begin{equation*}\label{J1x0}
\mathcal{M}(a,u_{01},u_{02}) 
=\dfrac{1}{2}\displaystyle \int_{t_1}^{t_2} |\eta(t) - \partial_x u^{a,u_{01},u_{02}}(1,t)|^2\, dt
   +
 \dfrac{1}{2}\displaystyle \int_{t_1}^{t_2} |\zeta(t) - \partial_x v^{a,u_{01},u_{02}}(1,t)|^2\, dt.
\end{equation*}
\noindent
\begin{test} \label{test5} 
We will take $\theta = 1.5$, $\alpha=1$, $\beta=1$, $T=2$, $t_1 = 0$, $t_2= T$, $u_{01ini} = 0.5$, $u_{02ini}= 1.8$ and $aini=0.1$ as initial guesses to recover the desired value of $a_d=0.5$, $u_{01d}= 1$, $u_{02d}= 2$ using the minimization algorithm. 
The numerical results can be seen in Figures~\ref{fig.aTest5}, ~\ref{fig.aTest5BB} and~\ref{fig.fvalTest5}. In Figure~\ref{fig.aTest5}, the stars and the round points correspond to the iterations during the optimization algorithm in the computation of $u_{01}$ and $u_{02}$, respectively.
In Figure~\ref{fig.aTest5BB}, the round points represent the iterations during the optimization algorithm in the computation of $a$.
With the solid line in Figure~\ref{fig.fvalTest5}, we have represented the evolution of the cost during iterations.  The algorithm does not converge well and we can see that the value of $u_{01}$ is not recovering properly.  The value of the functional does not become small: this indicates that, in this case, we are not able to obtain a solution of the inverse problem.
However, for the value of $u_{02}$ defined on the side where we take the distributed measurements, we get a better approximation. This suggests that the inverse problem might be solved in the right interval $(a,1)$, as we expect from the uniqueness Theorems \ref{theorem-strongly-right} and \ref{theorem2-uniq-system-2}.
\end{test}

\begin{figure}[]
\begin{minipage}[t]{0.49\linewidth}
\noindent
\includegraphics[width=0.98\linewidth]{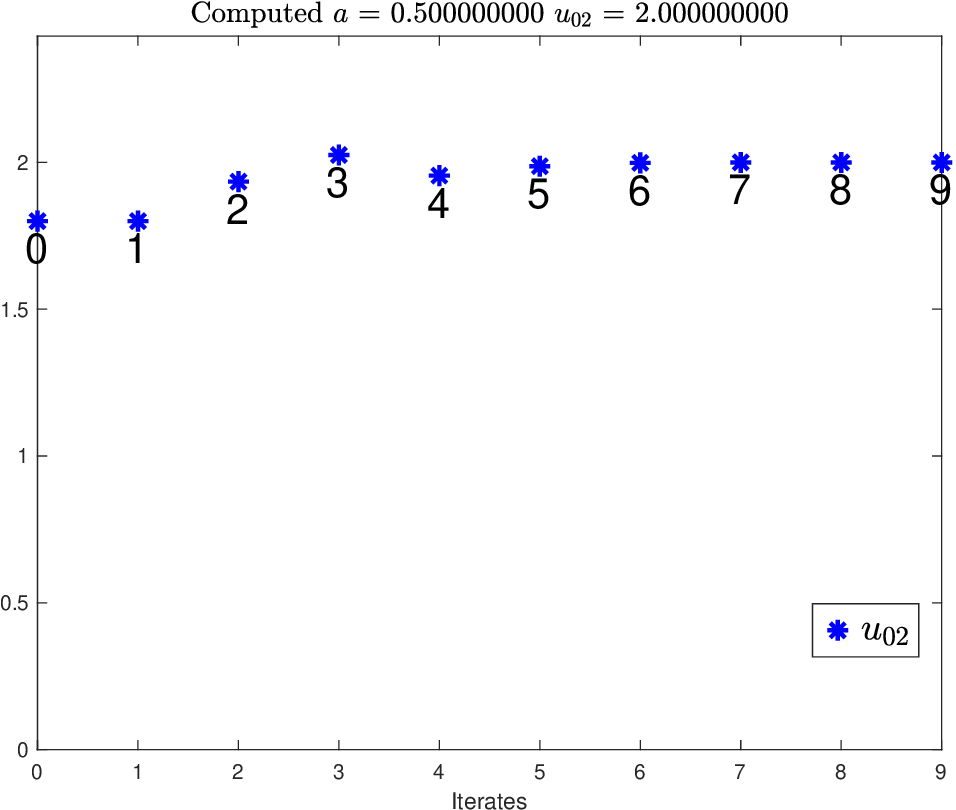}
\caption{Test~\ref{test6}, $\theta = 1.5$. Iterations in the computation of $u_{02}$ by \texttt{trust-region-reflective} algorithm, $aini=0.1$, $u_{02ini}= 1.8$. }
\label{fig.aTest6}
\end{minipage}
\hfill
\begin{minipage}[t]{0.49\linewidth}
\noindent
\includegraphics[width=0.98\linewidth]{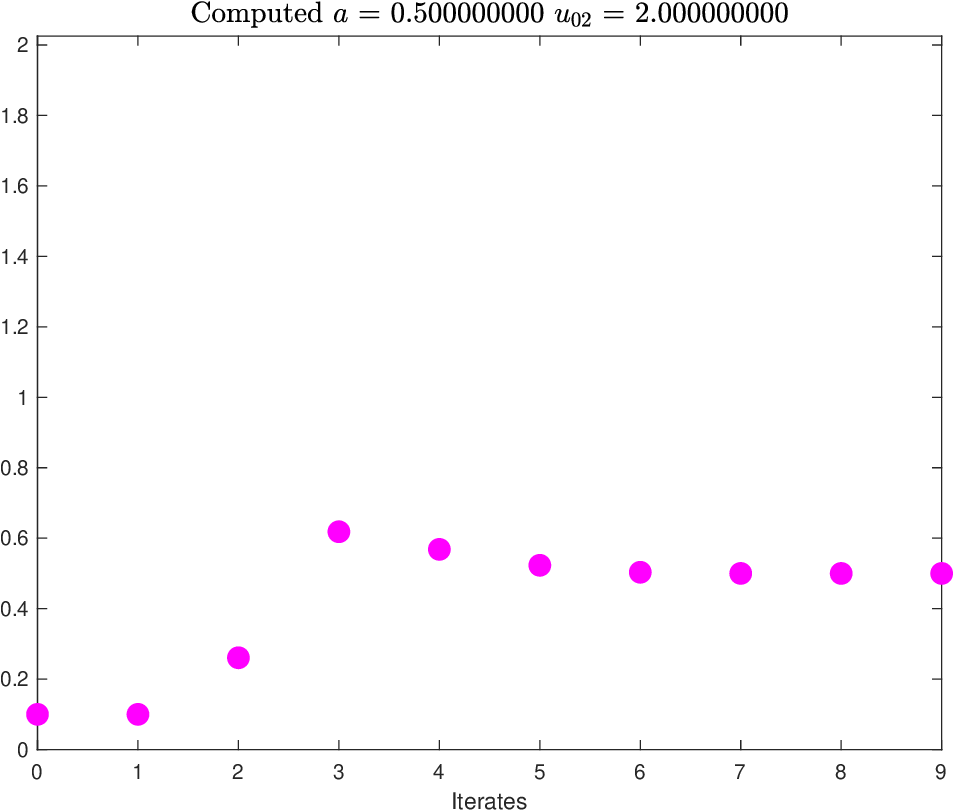}
\caption{Test~\ref{test6}, $\theta = 1.5$. Iterations in the computation of $a$ in \texttt{trust-region-reflective} algorithm, $aini=0.1$, $u_{02ini}= 1.8$.}
\label{fig.aTest6BB}
\end{minipage}
\end{figure}

\begin{figure}[h]
\centering
\includegraphics[width=0.49\linewidth]{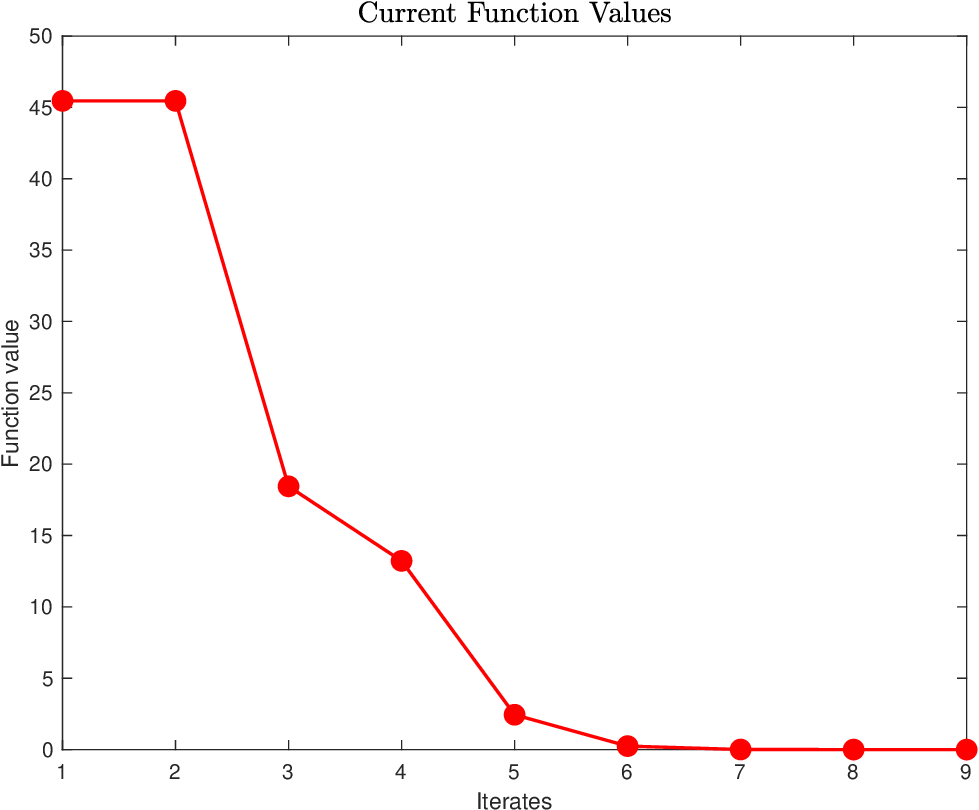}
\caption{Test~\ref{test6}, $\theta = 1.5$. Evolution of the cost in \texttt{trust-region-reflective} algorithm, $aini=0.1$ and $u_{02ini}= 1.8$.}
\label{fig.fvalTest6}
\end{figure}

In the next test, we will numerically obtain a result in line with the theoretical results \ref{theorem-strongly-right} and \ref{theorem2-uniq-system-2}.
Hence, we fix the initial data $u_{01}$ in $(0,a)$ and leave $a$ and the initial datum $u_{02}$ in $(a,1)$ as unknown, reconstructing them based on distributed measurements at $x=1$. 
The reformulation of the inverse problem is now as follows: 
\begin{equation*}
\left\{\begin{array}{l}
\text{Minimize  $\mathcal{K}(a,u_{02})$, }  \\[1mm]
\text{where $a\in \mathcal{U}^{a}_{ad} $ and $(u^{a,u_{02}},v^{a,u_{02}})$ satisfies \eqref{eq.cp},}
\end{array}\right.
\end{equation*}
\noindent
where $\mathcal{U}^{a}_{ad}$ is given by \eqref{uad} and $\mathcal{K} :(a,u_{02})\in \mathcal{U}^{a}_{ad}\times \mathds{R} \mapsto \mathds{R}$ is defined as follows:
\begin{equation*}\label{J1x0}
\mathcal{K}(a,u_{02}) 
=\dfrac{1}{2}\displaystyle \int_{t_1}^{t_2} |\eta(t) - \partial_x u^{a,u_{02}}(1,t)|^2\, dt
   +
 \dfrac{1}{2}\displaystyle \int_{t_1}^{t_2} |\zeta(t) - \partial_x v^{a,u_{02}}(1,t)|^2\, dt.
\end{equation*}
\begin{test} \label{test6} 
We will take $\theta = 1.5$, $\alpha=1$, $\beta=1$, $T=2$, $t_1 = 0$, $t_2= T$, $v_0=0$, $u_{01} = 1$, $u_{02ini}= 1.8$ and $aini=0.1$ as initial guesses to recover the desired value of $a_d=0.5$, $u_{02d}= 2$ using the minimization algorithm. 
The numerical results can be seen in Figures~\ref{fig.aTest6}, 	\ref{fig.aTest6BB} and~\ref{fig.fvalTest6}. In Figures~\ref{fig.aTest6} and \ref{fig.aTest6BB}, the stars and round points correspond to iterations in the computation of $u_{02}$ and $a$, respectively.  With the solid line in Figure~\ref{fig.fvalTest6}, we have represented the evolution of the cost.  The algorithm converges well, in particular, we can see an appropriate reconstruction of the value $u_{02}$ defined on the side where we take the distributed measurements.  We obtain the computed values  $a_c = 0.4999999996354925$, $u_{02c} = 2.0000000002802012$  and the cost $\mathcal{K}(a_c, u_{02c}) \approx 1.e-18$. 
\end{test}

\subsection{Degeneracy power and initial data reconstruction with distributed measurements}

\begin{figure}[h!]
\begin{minipage}[t]{0.49\linewidth}
\noindent
\includegraphics[width=0.98\linewidth]{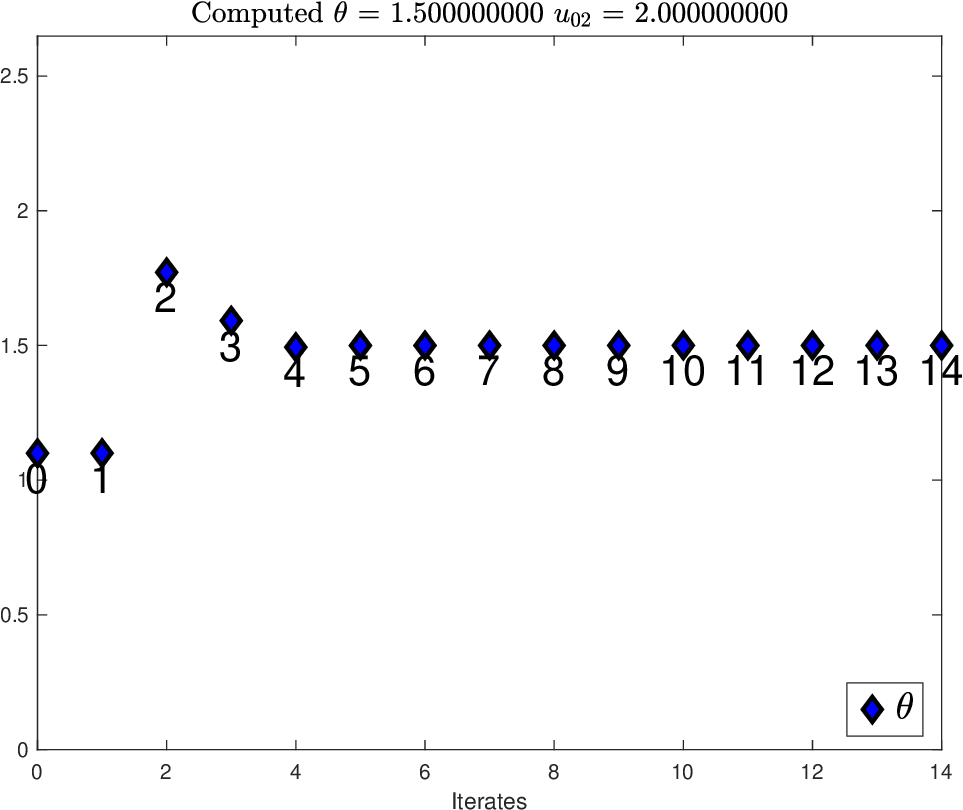}
\caption{Test~\ref{test6BIS}, $a = 0.5$. Iterations in the computation of $\theta$ by \texttt{trust-region-reflective} algorithm, $\theta_{ini}=1.1$, $u_{02ini}= 1.5$. }
\label{fig.thetaTest6BIS}
\end{minipage}
\hfill
\begin{minipage}[t]{0.49\linewidth}
\noindent
\includegraphics[width=0.98\linewidth]{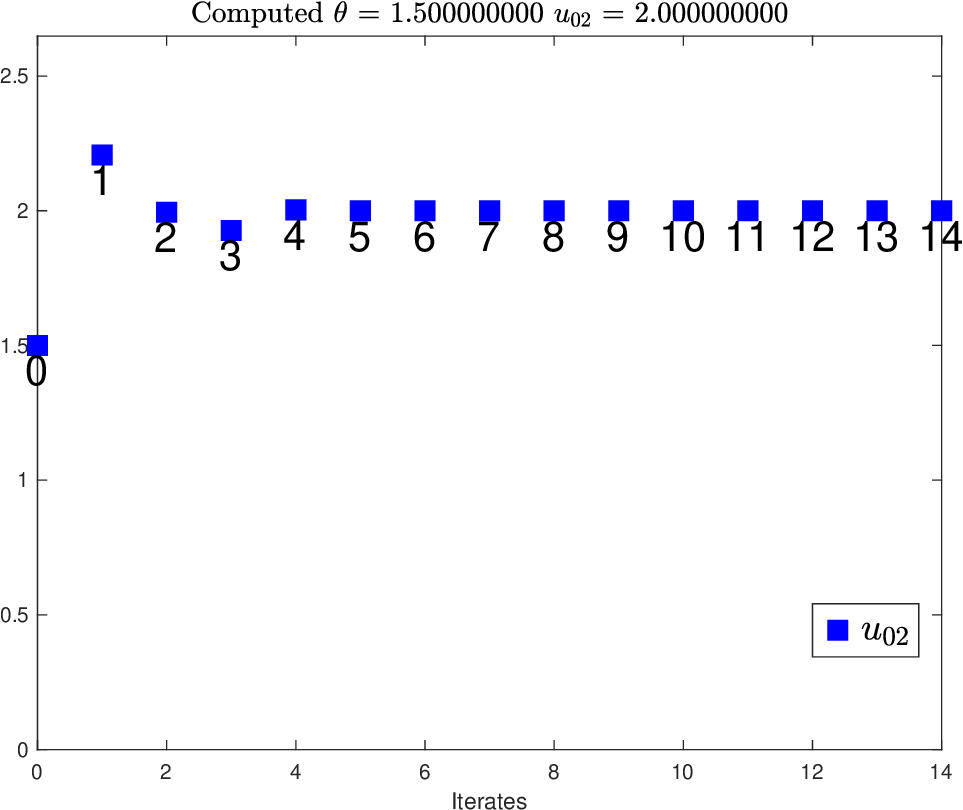}
\caption{Test~\ref{test6BIS}, $a = 0.5$. Iterations in the computation of $u_{02}$ by
\texttt{trust-region-reflective} algorithm, $\theta_{ini}=1.1$, $u_{02ini}= 1.5$.}
\label{fig.IniDataTest6BIS}
\end{minipage}
\end{figure}

\begin{figure}[h!]
\centering
\includegraphics[width=0.49\linewidth]{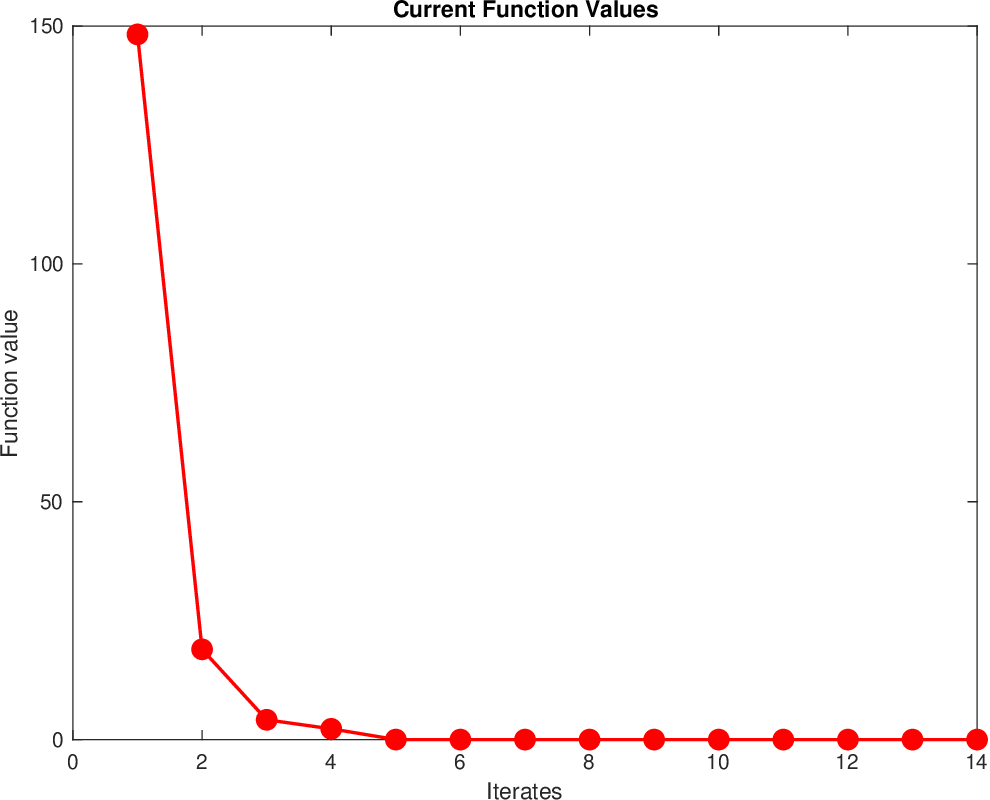}
\caption{Test~\ref{test6BIS}, $a = 0.5$. Evolution of the cost in \texttt{trust-region-reflective} algorithm, $\theta_{ini}=1.1$, $u_{02ini}= 1.5$.}
\label{fig.fvalTest6BIS}
\end{figure}

In this sub-section, we will give a numerical simulation in line with the results \ref{theorem-strongly-right-theta} and \ref{theorem2-uniq-system-theta}.
In particular, we assume that $u_0$ is of the form \eqref{inidata} and $v_0=0$, and fix the initial data $u_{01}=1$ in $(0,a)$, reconstructing the degeneracy power $\theta$ and the initial data $u_{02}$ in $(a,1)$ by distributed measurements at $x=1$.
Therefore, our goal is  to find $\theta\in [1,2)$ and initial data $u_{02}$ such that a solution to \eqref{eq.cp} satisfies
\begin{equation*}
\partial_x u(1,t) = \eta(t), 	\quad  \quad \partial_x v(1,t) = \zeta(t) \quad \text{ for } t \in (t_1,t_2), \quad 0\le t_1 \le t_2 \le T.
\end{equation*}
Now,  we reformulate the inverse problem as the following optimization problem: 
\begin{equation*}\label{optpb1x0BIS}
\left\{\begin{array}{l}
\text{Minimize  $\mathcal{N}(\theta,u_{02}) $, }   \\[1mm]
\text{where $\theta\in \mathcal{V}^{\theta}_{ad} $ and $(u^{\theta,u_{02}},v^{\theta,u_{02}})$ satisfies \eqref{eq.cp},}
\end{array}\right.
\end{equation*}
\noindent
where $\mathcal{V}^{\theta}_{ad} =\{ \theta: \theta\in [1, 2-\delta] \}$
and $\mathcal{N} : (\theta,u_{02})\in \mathcal{V}^{\theta}_{ad}\times \mathds{R} \mapsto \mathds{R}$ is defined as follows:
\begin{equation*}
\mathcal{N}(\theta,u_{02}) = \dfrac{1}{2}\displaystyle \int_{t_1}^{t_2} |\eta(t) - \partial_x u^{\theta,u_{02}}(1,t)|^2\, dt
+
 \dfrac{1}{2}\displaystyle \int_{t_1}^{t_2} |\zeta(t) - \partial_x v^{\theta,u_{02}}(1,t)|^2\, dt.
\end{equation*}
\noindent
\begin{test} \label{test6BIS}
We will take $a = 0.5$, $\alpha=1$, $\beta=1$, $T=4$, $t_1 = 0$, $t_2= T$, $u_{02ini}=1.5$ and $\theta_{ini}=1.1$ as initial guesses to recover the desired value of $\theta_d=1.5$, $u_{02d}=2$ by the minimization algorithm.  
The  numerical results can be seen in Figures~\ref{fig.thetaTest6BIS}, \ref{fig.IniDataTest6BIS} and~\ref{fig.fvalTest6BIS}. The rhombus-shaped points (Figure~\ref{fig.thetaTest6BIS}) and the squares (Figure~\ref{fig.IniDataTest6BIS}) correspond to iterations during the optimization algorithm in the computation of $\theta$ and $u_{02}$, respectively.  With the solid line in Figure~\ref{fig.fvalTest6BIS}, we have represented the evolution of the cost.  We obtain the computed value $\theta_c = 1.5000000000000084$, $u_{02c} = 2.0000000000000000$ and the cost $\mathcal{N}(\theta_c, u_{02c}) \approx 1.e-26$. 
\end{test}
\appendix
\section{Proof of Lemma \ref{lemma.Bessel}}

The proof of properties a), b) can be found in \cite{Watson}. For property c) see \cite{Abramowitz}.
With regard to property d), we have
\begin{equation*}
\begin{aligned}
    &\displaystyle\int_0^{j_{\nu,n}} s^{\nu+1} J_\nu(s) \, ds =\left[s^{\nu+1} J_{\nu+1}(s) \right]_0^{j_{\nu,n}}= j_{\nu,n}^{\nu+1} J_{\nu+1}(j_{\nu,n}) \\ &=-j_{\nu,n}^{\nu+1} J_\nu'(j_{\nu,n})+\dfrac{\nu j_{\nu,n}^{\nu+1}}{j_{\nu,n}}J_{\nu}(j_{\nu,n})= -j_{\nu,n}^{\nu+1} J_\nu'(j_{\nu,n}),
    \end{aligned}
\end{equation*}
where we have exploited properties a) and b).
The bounds on zeros of Bessel functions in f) and g) are given in \cite{Lorch-etal-08}.

\bigskip 

\noindent
\textbf{Acknowledgments:}
P.~C. and V.~D. were partially supported by the MUR Excellence Department Project awarded to the Department of Mathematics, University of Rome Tor Vergata, CUP E83C23000330006, by the PRIN 2022 PNRR-Project P20225SP98 “Some mathematical approaches to climate change and its impacts” (funded by the European Community-Next Generation EU, CUP E53D2301791 0001), and by Progetto di Ricerca di Ateneo 2024 Tor Vergata \textquotedblleft Control and Optimization for Multilayered Artificial
Neural Networks (COMANN)\textquotedblright, CUP E83C25000610005.
V.~D. was also supported by the PRIN Project 2022FPZEES “Stability in Hamiltonian Dynamics and
Beyond”.
P.~C. was also supported by the INdAM (Istituto Nazionale di Alta Matematica) Group for Mathematical Analysis, Probability and Applications. 
A.~D. was partially supported by Grants~PID2020-114976GB-I00 and~PID2024-158206NB-I00, funded by MICIU/AEI (Spain).
The authors warmly thank P. Martinez and M. Yamamoto for enlightening discussions and their precious comments.

\end{document}